\documentclass[
final,
11pt, a4paper, USenglish
]{article}

\providecommand\masterWhere{.}

\def\rootWhere{\masterWhere}

\def\arXivLibrary{\rootWhere/arXiv-library}

\usepackage{import}

\usepackage[USenglish]{babel}
\usepackage{import}

\usepackage[utf8]{inputenc}

\usepackage[T1]{fontenc}

\usepackage{mathscinet}

\usepackage{amsmath, amscd, amsthm}
\usepackage{thmtools}
\usepackage{amssymb}
\usepackage{mathtools}

\usepackage[shortcuts]{extdash}

\usepackage{tikz-cd}
\usepackage{tikz}

\usepackage{ifthen}

\newcommand{\ifisemptythenelse}[3]{%
  \if\relax\detokenize{#1}\relax
    #2%
  \else
    #3%
  \fi
}
\newcommand{\maybebrackets}[1]{%
				\ifisemptythenelse{#1}{}{{(#1)}}
				}

\usepackage{dsfont}
\usepackage[mathcal]{eucal}%

\newcommand\Eucal{\mathcal}%
\newcommand\ol{\overline}

\newcommand{\mathblackboardboldnew}{\mathbb}

\newcommand{\BN}{\mathblackboardboldnew N}

\newcommand{\BZ}{\mathblackboardboldnew Z}

\newcommand{\newmathcal}{\CMcal}

\newcommand{\calD}{\newmathcal{D}}

\newcommand{\calM}{\newmathcal{M}}

\newcommand{\rmB}{\mathrm{B}}

\newcommand{\rmE}{\mathrm{E}}

\newcommand{\rmK}{\mathrm{K}}

\newcommand{\rmS}{\mathrm{S}}

\newcommand{\frakC}{\mathfrak{C}}

\newcommand{\frakS}{\mathfrak{S}}

\newcommand{\qquote}[1]{``#1''}%
\newcommand{\roughly}[1]{``#1''}%

\newcommand{\introduce}{\textbf}%
\newcommand{\latin}{\emph}%
\newcommand{\buzzword}{\emph}%

\makeatletter
\newcommand\defpureword[2]{%
  \expandafter\newcommand\csname #1\endcsname[1]{\@pureword{#2}{##1}{#1}}
}
\newcommand{\@pureword}[3]{%
  \ifisemptythenelse%
  {#2}%
  {#1}%
  {
    \GenericError{}
    {Argument of pureword command not empty}
    {Use \ #3{} instead of \ #3.}
  }%
}
\newcommand\pureword[2]{\@pureword{#1}{#2}{command}}%
\makeatother

\newcommand{\Ealg}[1]{\rmE_{#1}}

\newcommand{\Eoperad}[1]{$\Ealg{#1}$-operad}

\defpureword{CWcomplex}{CW complex}

\defpureword{Ktheory}{$\rmK$-theory}

\defpureword{Sconstruction}{$\rmS$-construction}%

\defpureword{ConnesCC}{Connes' cyclic category}

\defpureword{inftycat}{$\infty$\=/category}
\defpureword{inftycats}{$\infty$\=/categories}
\defpureword{inftycategorical}{$\infty$\=/categorical}
\defpureword{pinftycat}{($\infty$\=/)category}
\defpureword{pinftycats}{($\infty$\=/)categories}
\defpureword{inftygrpd}{$\infty$\=/groupoid}
\defpureword{inftyoperad}{$\infty$\=/operad}
\defpureword{inftyoperads}{$\infty$\=/operads}
\defpureword{pinftyoperad}{{($\infty$\=/)}operad}
\defpureword{pinftyoperads}{{($\infty$\=/)}operads}

\newcommand{\nmcats}[2]{$(#1,#2)$\=/categories}

\newcommand{\resp}{{resp.\ }}
\newcommand\auxiliarytoggleablelatin\latin

\newcommand{\ie}{\auxiliarytoggleablelatin{i.e.}, }

\newcommand{\aka}{a.k.a.\ }

\newcommand{\non}[1]{non-{#1}}

\newcommand{\sub}[1]{sub{#1}}
\newcommand{\ary}[1]{$#1$\=/ary}
\newcommand{\twoptions}[2]{#1/#2}%

\newcommand{\twonames}[2]{{#1} and {#2}}
\newcommand{\threenames}[3]{{#1}, {#2} and {#3}}
\newcommand{\fivenames}[5]{{#1}, {#2}, {#3}, {#4} and {#5}}

\defpureword{Cech}{\v{C}ech}

\defpureword{TDyckerhoff}{T.~Dyckerhoff}
\defpureword{GJasso}{G.~Jasso}

\defpureword{htpyterminal}{homotopy terminal}
\defpureword{htpyinitial}{homotopy initial}
\defpureword{htpyinitialterminal}{homotopy \twoptions{initial}{terminal}}
\defpureword{wcontractible}{weakly contractible}
\defpureword{undercategory}{under-category}
\defpureword{overcategory}{over-category}
\defpureword{overundercategories}{\twoptions{over}{under}-categories}
\defpureword{homset}{Hom-set}
\defpureword{homfunctor}{Hom-functor}

\newtheoremstyle{newtheorem}{}{}{}{}{\bfseries}

\DeclareSymbolFont{extraup}{U}{zavm}{m}{n}
\DeclareMathSymbol{\varheart}{\mathalpha}{extraup}{86}
\DeclareMathSymbol{\vardiamond}{\mathalpha}{extraup}{87}
\DeclareMathSymbol{\varclub}{\mathalpha}{extraup}{84}
\DeclareMathSymbol{\varspade}{\mathalpha}{extraup}{85}

\newcommand{\somberend}{\ensuremath{\Diamond}}

\ifx\fancyends\undefined%
\newcommand{\theoremend}{\ensuremath{\square}}
\newcommand{\definitionend}{\somberend}
\newcommand{\exampleend}{\somberend}
\newcommand{\questionend}{\somberend}
\newcommand{\proofend}{\ensuremath{\blacksquare}}
\else
\newcommand{\theoremend}{\ensuremath{\varheart}}
\newcommand{\definitionend}{\ensuremath{\clubsuit}}
\newcommand{\exampleend}{\ensuremath{\diamondsuit}}
\newcommand{\questionend}{\ensuremath{\spadesuit}}
\newcommand{\proofend}{\ensuremath{\square}}
\fi

\theoremstyle{definition}

\newtheorem{Masterthm}{Masterthm}[subsection]

\declaretheorem[name=Construction ,style=definition,qed={\definitionend},sibling=Masterthm]{Cstr}
\declaretheorem[name=Construction ,style=definition,qed={\definitionend},unnumbered]{Cstr*}
\declaretheorem[name=Corollary,style=definition,qed={\proofend},sibling=Masterthm]{dCor}
\declaretheorem[name=Corollary,style=definition,qed={\theoremend},sibling=Masterthm]{Cor}
\declaretheorem[name=Definition,style=definition,qed={\definitionend},sibling=Masterthm]{Def}
\declaretheorem[name=Definition,style=definition,qed={\definitionend},unnumbered]{Def*}

\declaretheorem[name=Goal,style=definition,qed={\questionend},unnumbered]{Goal*}

\declaretheorem[name=Lemma,style=definition,qed={\theoremend},sibling=Masterthm]{Lem}

\declaretheorem[name=Proof,style=definition,qed={\proofend}, numbered=no]{Prf}

\declaretheorem[name=Proposition,style=definition,qed={\theoremend},sibling=Masterthm]{Prop}

\declaretheorem[name=Theorem,style=definition,qed={\theoremend},sibling=Masterthm]{Thm}
\declaretheorem[name=Theorem,style=definition,qed={\proofend},sibling=Masterthm]{dThm}
\declaretheorem[name=Theorem,style=definition,qed={\theoremend},unnumbered]{Thm*}

\theoremstyle{remark}
\declaretheorem[name=Remark ,style=remark,qed={\exampleend},sibling=Masterthm]{Rem}

\declaretheorem[name=Fact,style=remark,qed={\exampleend},unnumbered]{Fact*}

\declaretheorem[name=Example ,style=remark,qed={\exampleend},sibling=Masterthm]{Expl}


\numberwithin{equation}{section}

\declaretheorem[name=Theorem, style=definition, qed={\somberend}]{Theorem}
\declaretheorem[name=Theorem, style=definition, qed={\somberend},unnumbered]{Theorem*}

\declaretheorem[name=Corollary, style=definition, qed={\somberend},unnumbered]{Corollary*}

\declaretheorem[name=Observation, style=definition, qed={\somberend}, unnumbered]{Observation*}
\declaretheorem[name=Conjecture, style=definition, qed={\somberend}]{Conjecture}
\declaretheorem[name=Conjecture, style=definition, qed={\somberend},unnumbered]{Conjecture*}

\declaretheorem[name=Definition, style=definition, qed={\somberend}, unnumbered]{Definition*}

\declaretheorem[name=Remark ,style=remark,qed={\somberend},unnumbered]{Remark*}

\ifx\nopagestyle\undefined
\makeatletter
\newcommand\MyAtEndDocument[1]{%
  \g@addto@macro\My@EndDocument{#1}
}
\gdef\My@EndDocument{}
\AtEndDocument{\My@EndDocument}
\makeatother
\ifx\arXivLibrary\undefined %
\subimport{}{biblatex-style}

\else %

\date{}%
\MyAtEndDocument{
  \gdef\currenttitle{\bibname}
  \bibliography{\arXivLibrary}
}

\fi

\usepackage{geometry}
\usepackage{enumitem}%

\setlist{itemsep=-3pt, topsep=1pt}
\setlist[enumerate, 1]{label=(\arabic*), ref=(\arabic*)}

\usepackage{footnote}

\def\defaultdepthtoc{3}

\usepackage[titletoc]{appendix}				%
\usepackage[nottoc,notlot,notlof]{tocbibind} 	%
\addto\captionsUSenglish{%
  \renewcommand{\contentsname}%
    {Contents:}%
}
\addtocontents{toc}{\protect\setcounter{tocdepth}{\defaultdepthtoc}}
\newcommand\stoptoc{\addtocontents{toc}{\protect\setcounter{tocdepth}{-1}}}
\newcommand\resumetoc{\addtocontents{toc}{\protect\setcounter{tocdepth}{\defaultdepthtoc}}}

\makeatletter
\providecommand*\mytitle[1]{
  \gdef\@mytitle{\ifoptionfinal{#1}{#1 (work in progress)}}%
}
\providecommand*\myauthor[1]{\gdef\@myauthor{#1}}
\providecommand*\myaffiliation[1]{\gdef\@myaffiliation{#1}}
\providecommand*\myurl[1]{\gdef\@myurl{#1}}

\newcommand*\msnpError[1]{
  \GenericError{}
  {no my#1 provided}
  {Please call my#1 with an argument; do not try to def my#1}
}

\mytitle{\msnpError{title}}
\myauthor{\msnpError{author}}
\myaffiliation{\msnpError{affiliation}}
\myurl{\msnpError{url}}

\newcommand*\makemytitle{
  \title{\@mytitle}
  \author{\@myauthor\thanks{\@myaffiliation, \url{\@myurl}}}
  \maketitle
}

\gdef\currenttitle{\@mytitle}%

\usepackage{fancyhdr}
\usepackage{lastpage}

\fancypagestyle{tashi-style}{
  \fancyhf{}
  \fancyhead[LO]{\thepage/\pageref{LastPage}}
  \fancyhead[LE]{\@myauthor}
  \fancyhead[RO]{\currenttitle}
  \fancyhead[RE]{\thepage/\pageref{LastPage}}
  \fancyfoot[C]{}
}

\makeatletter

\fancypagestyle{plain}{%
  \fancyhf{}

}
\fi

\newcommand{\hra}{\hookrightarrow}

\newcommand{\lra}{\longrightarrow}
\newcommand{\lla}{\longleftarrow}

\newcommand{\lthra}%
{%
  \mathrel
  {
    \mathchoice
    {\mathrlap{\longrightarrow}\mkern 9mu\rightarrow}
    {\mathrlap{\longrightarrow}\mkern 9mu\rightarrow}
    {\mathrlap{\longrightarrow}\mkern 15mu\rightarrow}
    {\mathrlap{\longrightarrow}\mkern 15mu\rightarrow}
  }
}

\makeatletter

\newlength{\@minLengthArrowLengthOne}
\newlength{\@minLengthArrowLengthTwo}
\newlength{\@subscriptLengthForLongArrow}
\newlength{\@calculatedArrowLength}

\newcommand\@minLengthArrow [4][]{
  {
    \setlength{\@subscriptLengthForLongArrow}{#4}%
    \settowidth{\@minLengthArrowLengthOne}{\scriptsize$#1$}
    \settowidth{\@minLengthArrowLengthTwo}{\scriptsize$#2$}
    \pgfmathsetlength{\@calculatedArrowLength}{
      max
      ( \@minLengthArrowLengthOne
      , \@minLengthArrowLengthTwo
      , \@subscriptLengthForLongArrow
      )
    }
    \mathrel{
      #3
    [{\mathmakebox[\@calculatedArrowLength]{#1}}]
    {\mathmakebox[\@calculatedArrowLength]{#2}}
}  }
}

\newcommand\convertToMinLengthArrow[3][11pt]{
  \expandafter\newcommand\csname #2\endcsname[2][]{%
    \mathchoice
    {%
      \mathrel{\@minLengthArrow[##1]{##2}{#3}{#1}}
    }
    {%
      \mathrel{#3[##1]{##2}}
    }
    {%
      \mathrel{#3[##1]{##2}}
    }
    {%
      \mathrel{#3[##1]{##2}}
    }
  }
}

\makeatother

\convertToMinLengthArrow{xRa}{\xRightarrow}
\convertToMinLengthArrow{xLa}{\xLeftarrow}
\convertToMinLengthArrow{xLRa}{\xLeftrightarrow}
\convertToMinLengthArrow{xhla}{\xhookleftarrow}
\convertToMinLengthArrow{xhra}{\xhookrightarrow}
\convertToMinLengthArrow{xra}{\xrightarrow}
\convertToMinLengthArrow{xla}{\xleftarrow}
\convertToMinLengthArrow[13pt]{xlra}{\xleftrightarrow}
\convertToMinLengthArrow{xrat}{\xrightarrowtail}

\makeatletter
\newbox\xrat@below
\newbox\xrat@above
\newcommand{\xrightarrowtail}[2][]{%
  \setbox\xrat@below=\hbox{\ensuremath{\scriptstyle #1}}%
  \setbox\xrat@above=\hbox{\ensuremath{\scriptstyle #2}}%
  \pgfmathsetlengthmacro{\xrat@len}{max(\wd\xrat@below,\wd\xrat@above)+.6em}%
  \mathrel{\tikz [>->,baseline=-.75ex]
                 \draw (0,0) -- node[below=-2pt] {\box\xrat@below}
                                node[above=-2pt] {\box\xrat@above}
                       (\xrat@len,0) ;}}
\makeatother

\newcommand{\lrlas}{%
  \mathrel{\substack{\lra \\[-.65ex] \lla}}
}

\newcommand\adjarrows{\lrlas}%
\newcommand\ladjarrows\adjarrows%
\newcommand\radjarrows\adjarrows%

\newcommand\doublear[1]{\ar[#1, shift right]\ar[from=#1, shift right]}%

\usepackage{xstring}

\newcommand\intxt[1]{\qquad\text{#1}\qquad}

\newcommand\cdadjunctionOpt[5]{
  \IfEqCase{#1}{
    {u}{
      \ar[#1,shift left=2.5,"#2",/utils/exec={\pgfkeysalso{#4}}]
      \ar[#1,phantom,description,"\tiny{\dashv}"]
      \ar[from=#1,shift left=2.5,"#3",/utils/exec={\pgfkeysalso{#5}}]
    }
    {d}{
      \ar[#1,shift right=2.5,"#2"',/utils/exec={\pgfkeysalso{#4}}]
      \ar[#1,phantom,description,"\tiny{\dashv}"]
      \ar[from=#1,shift right=2.5,"#3"',/utils/exec={\pgfkeysalso{#5}}]
    }
    {r}{
      \ar[#1,shift left=2.5,"#2",/utils/exec={\pgfkeysalso{#4}}]
      \ar[#1,phantom,description,"\tiny{\bot}"]
      \ar[from=#1,shift left=2.5,"#3",/utils/exec={\pgfkeysalso{#5}}]
    }
    {l}{
      \ar[#1,shift right=2.5,"#2"',/utils/exec={\pgfkeysalso{#4}}]
      \ar[#1,phantom,description,"\tiny{\bot}"]
      \ar[from=#1,shift right=2.5,"#3"',/utils/exec={\pgfkeysalso{#5}}]
    }
  }
}

\newcommand\isCartesian[1][dr]{\ar[{#1}, phantom, description, very near start,"\lrcorner"]}
\newcommand\iscoCartesian[1][dr]{\ar[{#1}, phantom, description, very near end,"\ulcorner"]}
\newcommand\isbiCartesian[1][dr]{\ar[{#1},phantom, description, "\square"]}

\newcommand\cdsquare[9][]{
  \cdsquareOpt[#1]{#2}{#3}{#4}{#5}{"#6"}{"#7"'}{"#8"}{"#9"'}
}

\newcommand\cdsquareNA[5][]{\cdsquare[#1]{#2}{#3}{#4}{#5}{}{}{}{}}

\newcommand\cdsquareOpt[9][]{%
  \def\obja{#2}%
  \def\objb{#3}%
  \def\objc{#4}%
  \def\objd{#5}%
  \def\mora{\pgfkeysalso{#6}}
  \def\morb{\pgfkeysalso{#7}}
  \def\morc{\pgfkeysalso{#8}}
  \def\mord{\pgfkeysalso{#9}}
  \begin{tikzcd}[ampersand replacement=\&]
    \obja\ar[r,/utils/exec=\mora]\ar[d,/utils/exec=\morb]%
    \IfEqCase{#1}{%
      {C}{\isCartesian}%
      {cC}{\iscoCartesian}%
      {bC}{\isbiCartesian}%
    }
    \&%
    \objb\ar[d,/utils/exec=\morc]%
    \\%
    \objc\ar[r,/utils/exec=\mord]%
    \&%
    \objd%
  \end{tikzcd}%
}

\DeclareMathOperator{\Hom}{Hom}

\DeclareMathOperator{\Fun}{Fun}

\newcommand\Id{\mathrm{Id}}

\newcommand\sym{\mathrm{sym}}

\newcommand\blank{-}

\newcommand\inv{{-1}}%
\newcommand\op{\mathrm{op}}%

\newcommand\actson{\curvearrowright}%
\newcommand\basepoint\star%
\newcommand\disjunion{\mathbin{\dot\cup}}%
\newcommand\emptytuple\varnothing%
\newcommand\set[1]{\left\{#1\right\}}%
\renewcommand\emptyset\varnothing%

\newcommand\SG[1]{{\frakS_{#1}}}%
\newcommand\Sn{\SG n}%
\newcommand\cygrp[1]{\rquot \BZ {#1}}%
\newcommand\sphere[1]{S^{#1}}%
\newcommand\classB{\rmB}%

\newcommand\pushout[1]{\sqcup_{#1}}%
\newcommand\fiberproduct[1]{\times_{#1}}%
\newcommand\adjunit{\eta}%
\newcommand\ladjto\dashv%
\newcommand\radjto\vdash%
\newcommand\overcat[3][]{{{#2}_{/#3}}}%

\newcommand\grpdcore[1]{{#1}^{\simeq}}
\newcommand\rcone[1]{{#1}^\triangleright}%
 \newcommand\rquot[2]{
        \mathchoice
            {%
                \text{\raise1ex\hbox{${#1}$}}\Big/\lower1ex\hbox{${#2}$}%
            }
            {%
                #1\,/\,#2
            }
            {%
                #1/#2
            }
            {%
                #1/#2
            }
    }
\newcommand{\doubleslash}{\mathbin{
  \mathchoice{\Big/\mkern-10mu\Big/}%
    {/\mkern-6mu/}%
    {/\mkern-5mu/}%
    {/\mkern-5mu/}}}%

\newcommand{\qquot}[2]{
  \mathchoice
  {%
    \text{\raise1ex\hbox{${#1}$}}\doubleslash\lower1ex\hbox{${#2}$}%
  }
  {{#1}\doubleslash{#2}}
  {{#1}\doubleslash{#2}}
  {{#1}\doubleslash{#2}}
}

\newcommand\lquot[2]{
        \mathchoice
            {%
                \lower1ex\hbox{${#2}$}\Big\backslash\text{\raise1ex\hbox{${#1}$}}%
            }
            {%
                #2\backslash#1
            }
            {%
                #2\backslash#1
            }
            {%
                #1\backslash#2
            }
    }
    
 \newcommand\lrquot[3]{
        \mathchoice
            {%
                \lower1ex\hbox{${#2}$}\Big\backslash\text{\raise1ex\hbox{${#1}$}}\Big/\lower1ex\hbox{${#3}$}%
            }
            {%
                #2\backslash#1/#3
            }
            {%
                #2\backslash#1/#2
            }
            {%
                #1\backslash#2/#2
            }
    }

\newcommand*\noloc{%
        \nobreak
        \mskip6mu plus1mu
        \mathpunct{}%
        \nonscript
        \mkern-\thinmuskip
        {:}%
        \mskip2mu
        \relax
}

\newcommand\restr[3][]{{%
  \left.\kern-\nulldelimiterspace %
  {#2} %
  \vphantom{\big|} %
  \right|_{#3}^{#1} %
  }}

\newcommand{\category}[1]{\mathrm{\mathbf{#1}}}%
\newcommand{\cat}\category

\newcommand{\sSet}{\cat{sSet}}%
\newcommand\finset{\cat{Fin}}%

\newcommand\finsetneop{\finset_{\neq\emptyset}^\op}
\newcommand\finsetpop{\finset^\op_\star}

\newcommand{\Set}{\cat{Set}}%
\newcommand{\Cat}{\cat{Cat}}%
\newcommand\Spaces{\Eucal{S}}%

\newcommand\numD[1]{{[#1]}}%
\newcommand{\Dop}{{\Delta^{\op}}}

\newcommand\snerve[1]{\mathrm N_\Delta\maybebrackets{#1}}%
\newcommand\nerve[1]{\mathrm N\maybebrackets{#1}}%
\newcommand\catrealization[1]{\frakC\ifisemptythenelse{#1}{[\blank]}{[#1]}}%

\newcommand{\simplex}[1]{\Delta^{#1}}%
\defpureword{circlefunctor}{circle functor}%

\usepackage{hyperref}

\usepackage{ifdraft}
\usepackage[obeyFinal]{todonotes}
\AtEndDocument{\listoftodos}

\addto\extrasUSenglish{%
}

\mytitle{\twoSegal{} spaces as invertible \inftyoperad{}s}
\myauthor{Tashi Walde}
\myaffiliation{Rheinische Friedrich-Wilhelms-Universität Bonn}
\myurl{https://www.math.uni-bonn.de/people/twalde/}

\usepackage[all]{xy} %

\mathtoolsset{showonlyrefs}%
\usepackage{subfiles} %

\begin{document}
\makemytitle

\begin{abstract}
  
We exhibit the simplex category $\Delta$
and Segal's category $\Gamma$
as $\infty$-categorical localizations
of the dendroidal categories $\Omega_\pi$ and $\Omega$
introduced by Moerdijk and Weiss.
As an application we obtain an equivalence of $\infty$-categories
between invertible $\infty$-operads
and the $2$-Segal spaces
of Dyckerhoff and Kapranov.
Finally, we describe a cyclic version of the dendroidal category
and explain how it $\infty$-localizes to Connes's cyclic category $\Lambda$.

\end{abstract}

\tableofcontents

\section{Introduction}

\newcommand\exabcat{A}
\newcommand\Scstr{\rmS}
\newcommand\ScstrofA{\Scstr(\exabcat)}
\newcommand\ScAop[1]{\ScstrofA_{#1}}

\label{sec:introduction}
Higher category theory rests on the idea that one should replace
strictly associative composition of arrows
by composition laws which are only well-defined and associative up to
a coherent system of higher homotopies.
The importance of simplicial methods in higher category theory
stems mainly from the key fact that
the datum of an ordinary category can be faithfully repackaged
in a simplicial set, called its nerve.
It is well known that a simplicial set $\localX\colon \Dop\to \Set$ is
isomorphic to the nerve of a category if and only if
satisfies what are known as Rezk's \buzzword{Segal conditions}.
The category corresponding to $\localX$ has $\localX_{\numD 0}$
as its set of objects and $\localX_{\numD 1}$ as its set of morphisms;
composition of morphisms is defined by the span
\begin{equation}\label{span:segalcomposition}
  \mu\colon\localX_{\numD 1}\fiberproduct{\localX_{\numD 0}}\localX_{\numD 1}\xla{\cong}
  \localX_{\numD 2}\lra\localX_{\numD 1},
\end{equation}
where the left pointing map is guaranteed to be a bijection
by the first of the Segal conditions.
It is Rezk's fundamental insight~\cite{Rezk2001}
that one can model $(\infty,1)$-category as simplicial \emph{spaces}
which satisfy the correct homotopy coherent analog of the Segal conditions,
obtained by replacing bijections of sets by weak equivalences of spaces
and fiber products by their homotopy coherent counterparts;
the contractible (homotopy) fibers of the left pointing map
in~\eqref{span:segalcomposition} parameterize the choices of composition.

\nameDK{}~\cite{DyckerhoffKapranov2012} study the case where
the first map in the span~\eqref{span:segalcomposition}
is not an equivalence anymore.
In this case one can still interpret $\mu$ as
a \roughly{multi-valued composition law},
where the space of possible results of a composition
is parameterized by the possibly \non{contractible} or even empty fibers
of the first map in the span~\eqref{span:segalcomposition}.
This multi-valued composition law is unital and associative
(up to coherent homotopies)
precisely if the simplicial object $\localX$
satisfies the \buzzword{\twoSegal{} conditions}\footnote{
  \twoSegal{} spaces were also introduced independently
  by \nameGCKT{}~\cite{GCKT2018a,GCKT2018b,GCKT2018c}
  under the name \buzzword{decomposition spaces}.
}---a
weakening of Rezk's Segal conditions.

The main source of examples of \twoSegal{} spaces---apart
from all ordinary Segal spaces---is
\buzzword{Waldhausen's \Sconstruction{}}~\cite{Waldhausen1985},
which assigns to a suitable ($\infty$-)category $\localC$
a \twoSegal{} simplicial space $\wS\localC$
(see \autoref{expl:WaldhausenS}).
While Waldhausen was originally interested in the \emph{homotopical}
meaning of the \Sconstruction{}---the homotopy groups of $\wS\localC$
compute the algebraic \Ktheory{} of $\localC$---,
it turns out that the \Sconstruction{} also carries
interesting \emph{algebraic} information:
under suitable finiteness assumptions,
one can turn the simplicial space $\wS\localC$ into
the so called \buzzword{Hall algebra} of $\localC$
by an appropriate linearization procedure.
In this context, the \twoSegal{} property enjoyed by $\wS\localC$
can be seen to be directly responsible for the unitality and associativity
of the multiplication in the Hall algebra.
Variants of Hall algebras, such as the cohomological Hall algebra
of \twonames{Kontsevich}{Soibelmann}~\cite{KS11}
or the derived Hall algebra of Toën~\cite{Toen2006},
can be obtained by considering variants of this construction;
see \cite{Dyckerhoff2018} for a survey on this perspective.
\nameDK{} also recover classical convolution algebras
such as the \twonames{Iwahori}{Hecke} algebra
as linearizations of certain \twoSegal{} spaces.
Hall and Hecke algebras play an important role in representation theory,
for instance due to their close connection to quantum groups.

When constructing (strictly) associative algebras out of \twoSegal{} spaces,
one really only needs the $3$-skeleton of these simplicial spaces
and the corresponding truncated version of the \twoSegal{} conditions.
It is thus natural to ask:
What precisely is the higher algebraic structure
encoded in a \twoSegal{} space?
In this paper we establish the following theorem
(see \autoref{cor:invinfopsSegalspaces})
 which provides the first complete answer to this question.

\begin{Theorem}
  \label{thm:invertibleoperads}
  There is a canonical equivalence between
  \begin{itemize}
  \item the \inftycat{} of \twoSegal{} spaces and
  \item the \inftycat{} of \buzzword{invertible
    \inftyoperad{}s}\footnote{
     \label{fn:colored_nonsym}
      colored, \non{symmetric}
    }.\qedhere
  \end{itemize}
\end{Theorem}

The theory of \inftyoperad{}s,
originally introduced in the setting of
algebraic topology by May~\cite{May72} and Boardman--Vogt~\cite{BV73}
to study the algebraic structure of iterated loop spaces,
has since become a fundamental organizational tool
in the study of higher algebraic structures.
Roughly speaking,
an operad is a generalized category which admits
not just morphisms $x\to y$ from one object to another,
but also
\roughly{many-to-one} morphisms
$(x_1,\dots, x_n)\to y$,
called \buzzword{operations},
together with suitably associative composition laws
(see \autoref{def:coloredoperad}).

An operad is called \introduce{invertible}
(see \autoref{def:invertible-operad})
if each operation can uniquely be decomposed into other operations,
as long as the shape of this decomposition is specified in advance;
more precisely, we require that
each \ary{1} operation is the identity
and that,
after fixing
$0\leq i\leq j\leq n$,
each \ary{n} operation
$(x_1,\dots, x_n)\to z$
can be written uniquely
as a composition of two operations
$
(x_{i+1},\dots,x_{j})\to y
$
and
$
(x_1,\dots,x_{i},y,x_{j+1},\dots,x_n)\to z
$.
A trivial example of an invertible operad
is the commutative operad
which has a unique operation of each arity.
More interestingly,
there is, for each abelian category $\exabcat$,
an invertible operad $\ScstrofA$---%
corresponding to the aforementioned Waldhausen \Sconstruction{}
under the equivalence of \autoref{thm:invertibleoperads}---%
whose colors and \ary{1} operations are the objects of $\exabcat$
and whose \ary{2} operations are short exact sequences
(see \autoref{expl:WaldhausenS}).

The passage from operads to \inftyoperad{}s is analogous to the passage from
categories to \inftycats{} and arises by replacing strict composition of
operations by composition laws which are only well-defined and associative up to
a coherent system of higher homotopies.
To study \inftyoperad{}s we use the
convenient framework of dendroidal spaces introduced by
\nameMW{}~\cite{MW07} and later developed further by
\nameCM{}~\cite{CM11,CM13}.
In this framework the simplex category
$\Delta$ is replaced by a bigger category $\planarTrees$ of plane rooted trees
whose definition we recall in \autoref{sec:planarrootedTrees}.
Generalizing Rezk's ideas from
the simplicial case,
\nameCM{} observe that operads are identified
via a dendroidal version of the nerve functor
with \buzzword{dendroidal sets} $\planarTreesop\to\Set$
satisfying the dendroidal analog of the Segal conditions
(see~\autoref{def:SegalCM}).
More generally,
they show that \inftyoperad{}s are
modeled by (complete\footnote{
Completeness is an additional technical condition
which will be vacuous in the cases we consider.})
Segal dendroidal spaces.

The equivalence in~\autoref{thm:invertibleoperads}
is constructed by pulling back along an explicit functor
\begin{equation}
  \treerootpl\colon\planarTrees\lra\Delta
\end{equation}
(see \autoref{sec:treerootpl})
of ordinary categories,
which we prove to be an \inftycategorical{} localization in the following
sense:
\newcommand{\classofbp}{S}
There is an explicit class $\classofbp$ of
maps in $\planarTrees$ which are sent by $\treerootpl$ to equivalences in
$\Delta$ and,
moreover,
$\treerootpl$ is universal with this property among all
functors of \inftycats{}.
More precisely,
we have the following result
(see \autoref{thm:localizationW}).

\begin{Theorem}
\label{thm:intro:localization}
  Let $\localC$ be an \inftycat{}.
  The functor
  \begin{equation}
    \label{eq:intro-localiz}
    \treerootpl^\star\colon\Fun(\Delta,\localC)\lra\Fun(\planarTrees, \localC)
\end{equation}
  induced by $\treerootpl$ is fully faithful;
  the essential image is spanned by those functors
  $\planarTrees\to\localC$ which send all maps in $\classofbp$ to
  equivalences in $\localC$.
\end{Theorem}

\autoref{thm:invertibleoperads} follows from~\autoref{thm:intro:localization}
(after passing to opposite categories)
by observing that $\treerootpl^\star$
identifies \twoSegal{} simplicial objects in its domain
with (complete) Segal dendroidal objects in its essential image.\\

It is often worthwhile to enhance simplicial objects with
\roughly{additional symmetries}.
In this article we consider the following two main examples:

\begin{enumerate}
\item\label{it:enh_symmetric}
  Segal's \buzzword{special $\Gamma$-spaces}~\cite{Segal1974}---used
  to model the homotopy theory of connective spectra---can be seen
  as Segal simplicial spaces $\localX$ enhanced by compatible actions
  \begin{equation}\Sn\actson\localX_n\end{equation}
  of the symmetric groups.
\item\label{it:enh_cyclic}
  Cyclic symmetries on $\xdc$
  are encoded by lifts of $\localX$ to
  \buzzword{\ConnesCC{}} $\Lambda\supset\Delta$
  and described informally by a compatible system of actions
  \begin{equation}\mathrm{C}_{n+1}\actson\localX_n\end{equation}
  by cyclic groups.
  \twoSegal{} cyclic objects play a central role in Dyckerhoff--Kapanov's
  construction~\cite{DyckerhoffKapranov2018} of topological Fukaya categories of surfaces.
\end{enumerate}

One important feature of our proof of \autoref{thm:invertibleoperads}
is that it can be generalized to clarify how
cyclic (\resp symmetric) enhancements of \twoSegal{} spaces
correspond precisely to
cyclic (\resp symmetric) structures on
the corresponding invertible \inftyoperad{}s.
To do this we consider two variants of the category $\planarTrees$ of plane
rooted trees:
\begin{enumerate}
\item
  The category $\symTrees$ is precisely the category $\Omega$ of
  \nameMW{}. The objects of $\symTrees$ are rooted trees
  (without a chosen plane embedding);
  by the work of \nameCM{}~\cite{CM13},
  (complete) Segal presheaves on $\symTrees$
  are known to model \emph{symmetric} \inftyoperads{}.
\item
  By slightly modifying a construction of \twonames{Joyal}{Kock}~\cite{JK09},
  we introduce the category $\cycTrees$ of plane rootable trees
  (see \autoref{sec:cyclicTrees});
  it is expected\footnote{
    For instance, see \cite[Remark 6.9]{DH18} for a precise conjecture.
  }
  that (complete) Segal presheaves on $\cycTrees$
  are a model for \emph{cyclic} \inftyoperads{}.
\end{enumerate}
These categories of trees come equipped with canonical functors
\begin{equation*}
  \treerootsym\colon\symTrees\lra\Gamma
  \intxt{and}
  \treerootcyc\colon\cycTrees\lra\Lambda
\end{equation*}
(see \autoref{sec:symmetrictrees} and \autoref{sec:cyclicTrees}).
Our methods directly generalize to obtain the following version
of \autoref{thm:intro:localization} and \autoref{thm:invertibleoperads}
(see \autoref{thm:localizationWF},
\autoref{thm:localizationWL} and \autoref{rem:main-for-Gamma}).

\begin{Theorem}
  \label{thm:intro-Gamma-Lambda}
  The functors $\treerootsym$ and $\treerootcyc$
  are \inftycategorical{} localizations.
  Moreover the functor $\treerootsym$
  induces an equivalence of \inftycats{} between:
  \begin{itemize}
  \item
    \twoSegal{} $\Gamma$-spaces and
  \item
    invertible symmetric \inftyoperads{}.
    \qedhere
  \end{itemize}
\end{Theorem}

Since the localization functor $\treerootcyc$
identifies \twoSegal{} cyclic objects
with invertible Segal dendroidal objects,
\autoref{thm:intro-Gamma-Lambda} also implies the following conjecture
if we assume the conjectural existence of a
complete Segal cyclic dendroidal model for cyclic \inftyoperads{}
(see \autoref{rem:careful-with-Lambda}).

\begin{Conjecture}
  The functor $\treerootcyc$ induces an equivalence between
  \twoSegal{} cyclic spaces and
  invertible cyclic \inftyoperads{}.
\end{Conjecture}

\begin{Rem}
  The functor $\treerootsym\colon\symTrees\to\Gamma$
  was already considered by
  \nameBM{}~\cite[Theorem~1.1]{BM17};
  their main theorem states that this functor
  induces an equivalence between
  the \inftycat{} of special $\Gamma$-spaces and
  the \inftycat{} of what they call
  \buzzword{covariantly fibrant} complete Segal dendroidal spaces.
  We obtain their equivalence---%
  as well as the obvious variants for $\Lambda$ and $\Delta$---%
  by restricting our equivalences to the appropriate full subcategories
  (see \autoref{cor:Segalvscovariantfibration}).
\end{Rem}

\begin{Rem}
  Throughout this article
  we write \qquote{\twoSegal{}}
  to denote what \nameDK{} originally called \qquote{unital \twoSegal{}}.
  This is justified by the recent observation of
  \fivenames{Feller}{Garner}{Kock}{Proulx}{Weber}~\cite{FGKUPW2019}
  that unitality follows automatically from the \twoSegal{} conditions.
\end{Rem}

\begin{Rem}
  \autoref{thm:intro:localization} makes it possible to
  construct homotopy-coherent simplicial objects
  by specifying (possibly strict) dendroidal objects
  which send certain maps to weak equivalences.
  While this is easier \latin{a priori},
  the author does not know of any new simplicial objects that arise this way.
  When it comes to \twoSegal{} spaces, one should probably
  not expect new examples to arise from our result:
  first,
  because most operads appearing in the literature are not invertible
  and second, because every \twoSegal{} space
  can already be constructed by a generalized version
  of Waldhausen's \Sconstruction{}~\cite{BOORS2018}.
  Therefore, the results of this article should not be seen as a
  way to construct new \twoSegal{} spaces
  but rather as a new way of repackaging the higher algebraic structure
  encoded in such an object.
  This operadic perspective makes available tools and generalizations
  that were not evident in the original theory:
  While it is, for instance, not immediately obvious how to define
  \twoSegal{} objects with values
  in a general (not necessarily Cartesian) symmetric monoidal \pinftycat{},
  the definition of invertible \pinftyoperads{}
  directly generalizes to this setting;
  moreover, one can now hope to obtain new information about
  a \twoSegal{} space by studying
  algebras over the associated \inftyoperad{}.
\end{Rem}

\begin{Rem}
  Recently, a different algebraic interpretation of \twoSegal{} spaces
  was given by Stern~\cite{Stern2019}, who identified
  the \inftycat{} of \twoSegal{} objects in $\localC$
  with an \inftycat{} of algebras in correspondences in $\localC$.
  Similarly, Stern shows that \twoSegal{} cyclic objects
  are identified with Calabi-Yau algebras in correspondences.
\end{Rem}

\stoptoc

\subsection{Acknowledgments}

This work was done during my Ph.D.~studies at the Hausdorff Center for
Mathematics (HCM);
I am very grateful to my supervisors Catharina Stroppel and Tobias Dyckerhoff
for supporting and encouraging me during this time.
Furthermore, I would like to thank
Ieke Moerdijk and anonymous referees
for their very valuable feedback on previous versions of this paper.
This research was supported with a Hausdorff
scholarship by the Bonn International Graduate School for Mathematics (BIGS)
and funded by the Deutsche Forschungsgemeinschaft
(DFG, German Research Foundation)
under Germany's Excellence Strategy - GZ 2047/1, Projekt-ID 390685813.

\subsection{\inftycategorical{} conventions}

We use \inftycats{} (\aka quasi-categories)
as developed by Joyal~\cite{Joyal2002} and Lurie~\cite{Lurie2009}
as our preferred model for \nmcats{\infty}{1}.
Unless we are using specific technical results about simplicial sets,
we will, however, employ a rather high level language:
for instance, we treat each ordinary category $D$ as an \inftycat{}
by identifying it with its nerve $\nerve{D}$;
thus we usually write $D\to\localC$
for a functor to an \inftycat{} $\localC$,
and not $\nerve{D}\to\localC$.
We use the notation $\Fun(\localC', \localC)$ for the \inftycat{} of
functors $\localC'\to\localC$;
if $\localC'$ is (the nerve of) an ordinary category $D$
then we just write $\Fun(D,\localC)$.

\resumetoc

\section{The localization functors}

\label{sec:localizationfunctor}
Recall that the simplex category $\Delta$ is the
category of finite \non{empty} linearly ordered sets
and weakly monotone maps between them;
when convenient we identify $\Delta$ with its skeleton consisting of the
standard ordinals $\n=\set{0<\dots<n}$.

\subsection{The category $\planarTrees$ of plane rooted trees}

\label{sec:planarrootedTrees}
We recall some basic facts about (colored, \non{symmetric}) operads and the
category $\planarTrees$ of plane rooted trees as introduced by \nameMW{}~\cite{MW07}.

\begin{Def}\label{def:coloredoperad}
  A \introduce{colored, \non{symmetric} operad}
  (or \introduce{operad} for short)
  $\localO=(\localO, O, \ocomp)$ consists of
  \begin{itemize}
	\item a collection $O$ of \introduce{objects} (or \introduce{colors}),
	\item given colors $x_1, \dots, x_n,y\in O$, a set $\localO(x_1, \dots, x_n;
    y)$ of \introduce{$n$-ary operations} from $(x_1,\dots,x_n)$ to $y$ and
	\item for each $k, n_1, \dots, n_k \in\BN$ and colors $x^i_{j_i}, z \in O$
    (for $0 \leq j_i \leq n_i$, $0 \leq i \leq k$), a \introduce {composition
      map}
    \begin{align}
      \label{eq:operad-big-composition}
      \left(
      \coprod\limits_{y_1,\dots,y_k\in O}
      \big(
      \localO(x_1^1, \dots, x_{n_1}^1; y_1)
      \times\dots\times
      \localO(x_1^k, \dots, x_{n_k}^k; y_k)
      \big)
      \times
      \localO (y_1, \dots, y_k; z)
      \right)
      \\
      \xra{\ocomp}
      \localO(x_1^1,\dots, x^1_{n_1}, \dots, x_1^k, \dots, x_{n_k}^k; z)
    \end{align}
	\item
    a \introduce{unit map}
    \begin{equation}
      \label{eq:operad-big-unit}
      \ounit{}\colon O\lra\coprod_{x,y\in O}\localO(x;y)
    \end{equation}
    which assigns to each color
    $x\in O$
    the \ary{1} \introduce{identity operation}
    $\ounit{x}\in\localO(x;x)$
  \end{itemize}
  such that the obvious associativity and unitality conditions are satisfied.
  There is an obvious notion of a morphism of operads, we denote the resulting
  category of operads by $\Operads$.
\end{Def}

\begin{Rem}
  \label{rem:circle-i-compositions}
  By plugging suitable identity operations
  into the general composition law
  \eqref{eq:operad-big-composition}
  one can define the special compositions
  \begin{equation}
    \label{eq:circle-i-compositions}
    \circleis[j]{i}
    \colon
    \left(
    \coprod\limits_{y\in O}
    \localO(x_{i+1},\dots,x_{j};y)
    \times
    \localO(x_1,\dots,x_{i},y,x_{j+1},\dots,x_n;z)
    \right)
    \lra
    \localO(x_1,\dots,x_n;z)
  \end{equation}
  for all $0\leq i\leq j\leq n$
  and all $x_1,\dots,x_n,z\in O$.
  It is called $\circleis[j]{i}$ (the $j$ is left implicit)
  because the output of
  the first operation is inserted at position $i+1$ into the second.
  It is a easy to verify that every general composition map
  \eqref{eq:operad-big-composition}
  can be assembled as a suitable composition
  of such $\circleis[j]{i}$-compositions
  (for varying $i$ and $j$).
\end{Rem}

\begin{Rem}
  As originally introduced by Boardman--Vogt and May, an \qquote{operad} would
  be assumed to be mono-colored. Since there is no reason for us to single out
  this special case we will instead take operads to be colored by default.
  Moreover it is most convenient for us to reserve the word \qquote{operad} for
  the least structured situation and add further adjectives (e.g.\ symmetric or
  cyclic) whenever we equip our operad with extra structure (see also
  \autoref{sec:cyclicTrees} and \autoref{sec:symmetrictrees}). We warn
  the reader that this is a rather uncommon convention: most authors (including
  \nameMW{} and \nameCM{}) will define
  operads to be symmetric by default.
\end{Rem}

\begin{Rem}
  Each operad $(\localO, O,\ocomp)$
  has an underlying category with objects $x\in O$ and morphism sets
  $\localO(x;y)$. Conversely, each category can be viewed as an operad which has
  only $1$-ary operations. More precisely, we have an adjunction
  $\Cat\radjarrows\Operads$ with fully faithful left adjoint.
\end{Rem}

An object of $\planarTrees$ is called a \introduce{plane rooted tree} and
consist of a finite plane rooted trees in the usual graph-theoretic sense
together with a marking of some degree $1$ vertices including the root-vertex.
An edge between unmarked vertices is called internal, the other edges are called
external. The unique external edge connected to the root-vertex is called the
\introduce{root} (or output edge); an external edge attached to a marked
\non{root} vertex is called a \introduce{leaf} (or input edge).

\begin{Expl}
  \label{expl:some-trees}
  We depict some trees in $\planarTrees$,
  including the special tree $\degtree$, some corollas
  ($\corolla 0$, $\corolla 1$, $\corolla 3$)
  and two typical trees (of arity $3$ and $4$, respectively).
  
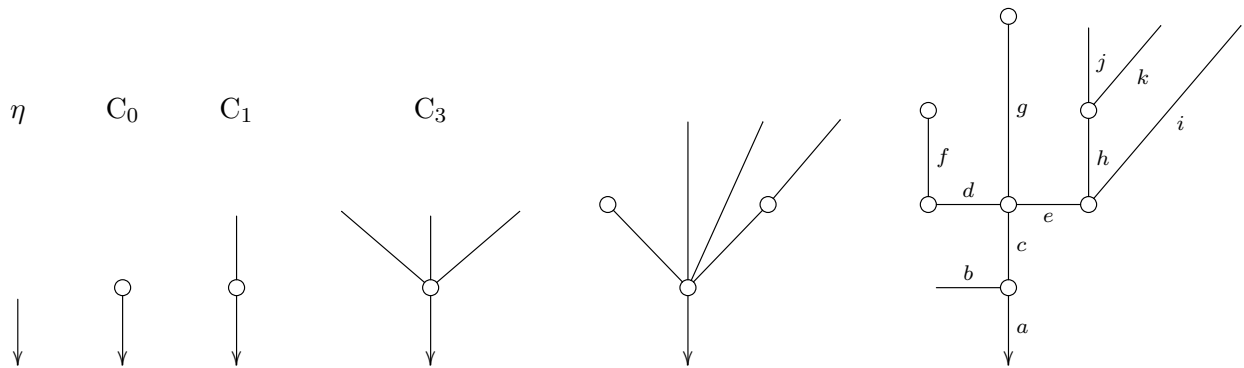
\begin{equation*}
    \xymatrix{
      &&&&&&																&&&&												&*+=+[o][F]{}\ar@{-}[dd]^{g}&\ar@{-}[d]^{j}&\ar@{-}[dl]^{k}&\ar@{-}[lldd]^{i}		&\\
      \degtree&\corolla 0&\corolla 1&&\corolla 3&&&																			\ar@{-}[dd]&\ar@{-}[ddl]&\ar@{-}[dl]&					*+=+[o][F]{}\ar@{-}[d]^{f}&&*+=+[o][F]{}\ar@{-}[d]^{h}&&\\
      &&\ar@{-}[d]&\ar@{-}[dr]&\ar@{-}[d]&\ar@{-}[dl]&							*+=+[o][F]{}\ar@{-}[dr]&&*+=+[o][F]{}\ar@{-}[dl]&&		*+=+[o][F]{}\ar@{-}[r]^{d}&*+=+[o][F]{}\ar@{-}[d]^{c}&*+=+[o][F]{}\ar@{-}[l]^{e}&&\\
      \ar@{->}[d]&*+=+[o][F]{}\ar@{->}[d]&*+=+[o][F]{}\ar@{->}[d]&&*+=+[o][F]{}\ar@{->}[d]&&		&*+=+[o][F]{}\ar@{->}[d]&&&							\ar@{-}[r]^{b}&*+=+[o][F]{}\ar@{->}[d]^{a}&&&\\
      &&&&&&																				&&&	&												&&&&\\
    }
  \end{equation*}
  The root is marked with a little arrow and drawn towards the bottom.
\end{Expl}

\begin{Rem}
  From now on we completely ignore the marked vertices of a tree and never speak
  of them again. Thus \qquote{vertex} always means \qquote{unmarked vertex}.
  When drawing trees, we omit the marked vertices and instead draw the external
  edges \roughly{towards infinity}.
\end{Rem}

The number of leaves of a tree is its \introduce{arity}. Each vertex of a tree
has some number (the \introduce{arity} of that vertex) of input edges and a
unique output edge (which is the one that points in the direction of the root).
The input edges of a vertex are linearly ordered left-to-right by the plane embedding.
We denote by $\degtree$ or ${\numD 0}$ the tree with only a single edge
(which is both the root and a leaf);
we denote by $\corolla \n$ or $\corolla n$ the
\introduce{$n$\=/corolla},
i.e.\ the unique $n$-ary tree with a single vertex.
Given two edges $e, e'$ in a plane tree $T$,
we say that $e$ is a \introduce{predecessor} of $e'$
and that $e'$ is a \introduce{successor} of $e$,
if the unique path in $T$ going from $e$ to the root of $T$ goes through $e'$;
note that every edge is a predecessor of the root.
Given two edges $d,e$ in $T$,
we say that $d$ \introduce{lies to the left} of $e$
and that $e$ \introduce{lies to the right} of $d$,
if there are successors $d'$ of $d$ and $e'$ of $e$
which are input edges at a common vertex $v$ and such that $e'$ lies (strictly)
to the left of $e$ with respect to the left-to-right linear order at $v$.
Observe that for any two edges $e,d$ we have the following two mutually
exclusive cases:
\begin{itemize}
\item
  $d$ is a successor or a predecessor of $e$ (this includes the case $d=e$) or
\item
  $d$ lies to the left or to the right of $e$.
\end{itemize}

\begin{Expl}
  In the last tree of \autoref{expl:some-trees}:
  The predecessors of the edge $e$ are $e$ itself, $h$, $j$, $k$ and $i$;
  the successors of $e$ are $e$ itself, $c$ and the root $a$.
  To the left of $e$ lie the edges $d$, $f$, $g$ and $b$;
  no edge lies to the right of $e$.
\end{Expl}

Each plane rooted tree $T$ gives rise to a free operad (also denoted by $T$):
it has a color for each edge of $T$ and
its operations are freely generated by the vertices of T
(an $n$-ary vertex is seen as an $n$-ary operation
from its input edges to its output edge).
A morphism in $\planarTrees$ between two trees is
defined to be a morphism of the corresponding operads.

\begin{Expl}\label{expl:somemapoftrees}
  Consider the following two plane rooted trees. The operad associated
  to the left tree has colors $\{a', a,c,d,e,f\}$ and three \non{unit} operations
  $s\colon a'\to a$ and $r\colon (e,f,c,d)\to a'$ and $r\ocomp
  s\colon(e,f,c,d)\to a$. The other one has colors $\{\ol a, \ol b,\ol c,\ol
  d,\ol e, \ol f,\ol g, \ol h\}$ and eleven \non{unit} operations ($t$, $u$, $v$,
  $w$ and all their composites).

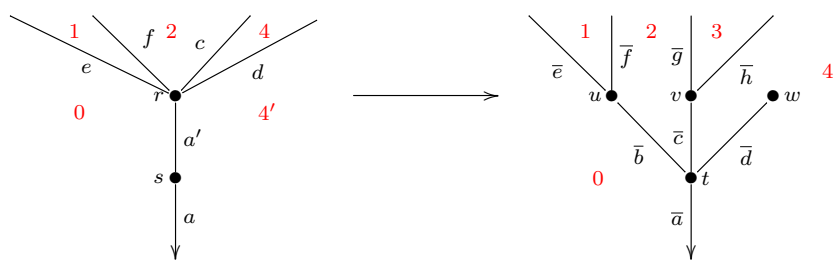
\begin{equation}\xymatrix{*{\,}\ar@{-}[rrd]_{e} \ar@{}[r]_{\,\,\,\,\,\,\,\,\,\,\,\,\color{red}{1}}\ar@{}[ddrr]_{\,\,\,\,\,\,\,\,\,\,\,\,\color{red}{0}}& *{\,}\ar@{}[rr]_{\color{red}{2}}\ar@{-}[rd]^{f} &  & *{\,}\ar@{-}[dl]_{c} & \ar@{-}[lld]^{d}\ar@{}[l]^{\color{red}{4}\,\,\,\,\,\,\,\,\,\,\,\,}\ar@{}[ddll]^{\color{red}{4'}\,\,\,\,\,\,\,\,\,\,\,\,}\\
 & \,\ar@{}[r]|{\,\,\,\,\,\,\,\,\,\,\,\,\, r} & *{\bullet}\ar@{-}[d]^{a'}\\
 & \,\ar@{}[r]|{\,\,\,\,\,\,\,\,\,\,\,\,\, s} & *{\bullet}\ar@{->}[d]^{a}\\
 &  & *{\,}}
\xymatrix{\\\,\ar[rr] & & \,} \xymatrix{
  *{\,}\ar@{-}[dr]_{\ol e}\ar@{}[r]_{\,\,\,\,\,\,\,\,\color{red}{1}} & *{\,}\ar@{-}[d]^{\ol f} \ar@{}[r]_{\color{red}{2}}&*{\,}\ar@{-}[d]_{\ol g} & *{\,}\ar@{-}[dl]^{\ol h}\ar@{}[l]^{\color{red}{3}\,\,\,\,\,\,\,\,}&\\
  \,\ar@{}[r]|{\,\,\,\,\,\,\,\,\,\,\,\,\, u} & *{\bullet}\ar@{-}[dr]_{\ol b}\ar@{}[dd]_{\color{red}{0}} \ar@{}[r]|{\,\,\,\,\,\,\,\,\,\,\,\,\, v}&*{\bullet}\ar@{-}[d]_{\ol c}\ar@{}[r]|{\,\,\,\,\,\,\,\,\,\,\,\,\,\,\,\,\,\,\,\,\,\,\,\,\,\,\,\,\,\,\,\, w}&   *{\bullet}\ar@{-}[dl]^{\ol d}\ar@{}[ru]_{\color{red}{4}}&\\
  &  & *{\bullet}\ar@{->}[d]_{\ol a} & \,\ar@{}[l]|{t\,\,\,\,\,\,\,\,\,\,\,\,\,\,}\\
  & & *{\,}} \end{equation} The depicted morphism is described on colors by $a'\mapsto \ol
a$, $a\mapsto \ol a$, $c\mapsto\ol c$ etc.\ and on generating operations by
$s\mapsto \ounit{\ol a}$ and $r\mapsto (u,\ounit{\ol c},\ounit{\ol d})\ocomp t$.
(The red numbers are for later reference.)
\end{Expl}

A \introduce{(planar) dendroidal object} in an \inftycat{} $\localC$ is
functor $\planarTreesop\to \localC$. We denote by
$\pldSet\coloneqq[\planarTreesop,\Set]$ the category of (planar) dendroidal
sets, i.e.\ dendroidal objects in $\Set$. Given a plane rooted tree $T$, we
denote by $\planarTrees[T]$ the dendroidal set represented by $T$. There is a
canonical fully faithful embedding $\Delta\hra\planarTrees$ of the simplex
category $\Delta$ by interpreting every linearly ordered set as a linear tree.
This embedding gives rise to an adjunction $\sSet\radjarrows \pldSet$ with fully
faithful left adjoint. The inclusion $\planarTrees\hra \Operads$ (which is full
by construction) gives rise to a realization/nerve adjunction
\begin{equation}\pldSet\adjarrows \Operads\noloc \dnerve{}\end{equation}
by the formula $\dnerve{\localO}\colon T\mapsto\Hom_{\Operads}(T, \localO)$,
which extends the usual adjunction
\begin{equation}\sSet\adjarrows\Cat\noloc \nerve{}.\end{equation}

\subsection{The localization functor $\treerootpl\colon\planarTrees\to\Delta$}

\label{sec:treerootpl}
Let us introduce the main player in our game.

\begin{Cstr}[Covariant description of $\treerootpl$]
  \label{cstr:treerootplcovariant}
  Each plane rooted tree $T\in \planarTrees$ (which we visualize with its
  external edges going towards infinity) partitions the plane into a set
  $\treerootpl T$ of \roughly{areas} which is linearly ordered clockwise
  starting from the root. It is straightforward to extend this assignment to a
  functor $\treerootpl\colon\planarTrees\to\Delta$.
\end{Cstr}

We give an alternative, more formal, construction of the functor $\treerootpl$
at the end of this section, see \autoref{cstr:treerootplcontravariant} below.

\begin{Expl}
  The functor $\treerootpl$ sends the morphism depicted in
  \autoref{expl:somemapoftrees} to the map $\{0,1,2,4,4'\}\to\{0,1,2,3,4\}$ in
  $\Delta$ which sends $i',i \mapsto i$.
\end{Expl}

\begin{Rem}
  Specifying two adjacent \roughly{areas} of a plane rooted tree
  $T\in\planarTrees$ uniquely determines an external edge of $T$ that separates
  them. If we write $\n\coloneqq \treerootpl T$ (where $n$ is the arity of $T$)
  then
  \begin{itemize}
  \item each minimal edge $\{i-1,i\}\hra \n$ (for $1\leq i\leq n$) corresponds
    precisely to a leaf of $T$ and
  \item the maximal edge $\{0,n\}\hra \n$ corresponds to the root of
    $T$.\qedhere
  \end{itemize}
\end{Rem}

\begin{Rem}
  Usually the category of trees is related to the simplex category by the
  inclusion $\Delta\hra\planarTrees$ of the linear trees. The composition
  $\Delta\hra \planarTrees\xra{\treerootpl}\Delta$ is constant with value ${\numD 1}\in
  \Delta$. The two occurrences of the category $\Delta$ in relation to the
  category $\planarTrees$ are in some sense \roughly{orthogonal}: the first is
  sensitive to the \roughly{height} of a tree, the second measures the
  \roughly{width}.
\end{Rem}

\begin{Def}
  A map of plane rooted trees is called \introduce{boundary preserving} if it
  maps the root to the root and each leaf to a leaf.
\end{Def}

\begin{Def}\label{def:invertibledendr}
  A \introduce{collapse map} in $\planarTrees$ is a boundary preserving map
  $\corolla\n\to T$ out of a corolla (where $n$ is the arity of $T$). A
  dendroidal object $\localX\colon \planarTreesop\to\localC$ in some
  \inftycat{} $\localC$ is called \introduce{invertible} if $\localX$
  maps all collapse maps to equivalences in $\localC$.
\end{Def}

\begin{Rem}
  A boundary preserving map $\alpha\colon T\to S$ of plane rooted trees induces
  a bijection between the leaves of $T$ and the leaves of $S$. Hence the functor
  $\treerootpl$ maps boundary preserving maps to isomorphisms.
\end{Rem}

\begin{Rem}
  The motivation for the word \qquote{invertible} in
  \autoref{def:invertibledendr} will become apparent in
  \autoref{sec:invertibleoperads} when we discuss invertible operads (in the
  sense of Dyckehoff and Kapranov~\cite{DyckerhoffKapranov2012})
  and show that an operad is invertible if and only if its nerve is an
  invertible dendroidal set (\autoref{lem:charinvertibleoperads}).
\end{Rem}

Here is one version of our main result which we explain and prove in
\autoref{sec:mainproof} below:

\begin{Thm}\label{thm:localizationWD}
  The functor $\treerootpl$ exhibits $\Delta$ as an \inftycategorical{}
  localization of $\planarTrees$ at the set of collapse maps.
\end{Thm}

Before going forward, we give a \roughly{contravariant} description of the
functor $\treerootpl$. This description is useful because unlike the covariant
one it can easily be adapted to the case of symmetric trees (see
\autoref{sec:symmetrictrees}).
Denote by $\boundedLinsets$ the following category:
objects are (possibly empty) linearly ordered sets;
a morphism $N\to M$ is a weakly monotone map
\begin{equation}\set{-\infty}\disjunion N\disjunion\set{+\infty}\to
\set{-\infty}\disjunion M\disjunion\set{+\infty}\end{equation}
which preserves $-\infty$ and $+\infty$
(where $-\infty$ and $+\infty$ are a new minimal and maximal element,
respectively).
It is an easy fact (sometimes known as Joyal duality)
that the category $\Delta$ is isomorphic to $\boundedLinsets^\op$
via the assignment (described here only on objects)
\begin{equation}\Delta\ni N\longmapsto
  \set{\text{\non{empty} proper initial segments of }N}\in \boundedLinsets^\op.\end{equation}

Using the identification $\Delta\simeq \boundedLinsets^\op$
we can give the following description of the functor
$\treerootpl\colon \planarTrees\to\boundedLinsets^\op$,
which is easily seen to be equivalent
to \autoref{cstr:treerootplcovariant}.

\begin{Cstr}[Contravariant description of $\treerootpl$]
  \label{cstr:treerootplcontravariant}
  To each plane rooted tree $T\in\planarTrees$ we associate the (possibly empty)
  linearly ordered set $\treerootpl T\in\boundedLinsets$ of its leaves. This
  association extends to maps in the
  following way:
  Given a map $\alpha\colon S\to T$ of trees,
  we need to define a map
  $\set{-\infty}\disjunion\treerootpl T\disjunion\set{+\infty}
  \to \set{-\infty}\disjunion\treerootpl S\disjunion\set{+\infty}$.
  We have no choice but to send
  $-\infty$ and $+\infty$ to $-\infty$ and $+\infty$, respectively.
  Denote by $\treerootof{S}$ the root of $S$ and let $a \in \treerootpl{T}$;
  there are three cases:
  \begin{itemize}
  \item
    If $a$ is a predecessor of $\alpha(\treerootof{S})$
    then there is a unique leaf $b$ of $S$
    such that $\alpha(b)$ is a successor of $a$;
    in this case we define
    $(\treerootpl\alpha)(a)\coloneqq b$ to be this unique leaf.
  \item
    If $a$ lies to the left of $\alpha(\treerootof{S})$
    then we define
    $(\treerootpl\alpha)(a)\coloneqq -\infty$.
  \item
    If $a$ lies to the right of $\alpha(\treerootof{S})$
    then we define
    $(\treerootpl\alpha)(a)\coloneqq +\infty$.
  \end{itemize}
  It is straightforward to verify that this assignment defines a
  functor $\treerootpl\colon\planarTrees \to\boundedLinsets^\op$.
\end{Cstr}

\begin{Expl}
  The map of trees from \autoref{expl:somemapoftrees} gets sent
  by $\treerootpl$ to the map
  \begin{equation}\set{-\infty,\ol e,\ol f,\ol g,\ol h,+\infty}
    \to\set{-\infty,e,f,c,d,+\infty}\end{equation}
  in $\boundedLinsets$
  given by $\ol e\mapsto e$, by $\ol f\mapsto f$ and by $\ol g,\ol h\mapsto c$.
\end{Expl}

\subsection{Symmetric operads and Segal's category $\Gamma$}

\label{sec:symmetrictrees}

Before moving on with the proof of our main localization theorem,
we briefly describe the analog construction
in the world of symmetric operads,
\ie operads equipped with compatible action of the symmetric groups
which interchange the input colors of an operation.

Denote by $\symTrees$ the category of symmetric rooted trees
(i.e.\ trees without a plane embedding),
defined as a suitable full subcategory of the category $\symOperads$
of symmetric operads;
this is the category of trees which
\nameMW{}~\cite[Section~3]{MW07} simply call $\Omega$.
Boundary preserving maps and collapse maps in $\symOperads$
are defined in the same way as for plane trees.

The symmetric analog of the simplex category
is Segal's category $\Gamma\coloneqq\finsetpop$,
the opposite of the category of finite pointed sets.
We define a functor $\treerootsym\colon\symTrees\to\Gamma$,
which is analogous to
$\treerootpl$ by adapting \autoref{cstr:treerootplcontravariant}:

\begin{Cstr}[The functor $\treerootsym$]
  We define the functor $\treerootsym\colon\symTrees\to\finsetpop=\Gamma$
  as follows:
  To each tree $T$ we
  assign the set of external edges which is pointed at the root.
  Given a morphism $\alpha\colon S\to T$ of rooted trees
  and a leaf $a$ of $T$ there is at most one leaf $b$ of $S$
  such that $\alpha(b)$ is a successor of $a$;
  we define
  $(\treerootsym\alpha)(a)\coloneqq b$ if such a $b$ exists and
  $(\treerootsym\alpha)(a)\coloneqq \basepoint$ otherwise.
\end{Cstr}

It is straightforward to show that $\treerootsym\colon\symTrees\to\finsetpop$ is
well defined and extends the functor $\treerootpl$ in the sense that the
following diagram commutes:
\begin{equation}\begin{tikzcd}
    \Operads\ar[d, "\sym"]\ar[r, hookleftarrow]&\planarTrees\ar[d]\ar[r, "\treerootpl"]& \Delta\ar[r,"{\simeq}", leftrightarrow]&\boundedLinsets^\op\ar[ld]\\
    \symOperads\ar[r, hookleftarrow]&\symTrees\ar[r,"\treerootsym"]&\finsetpop&
  \end{tikzcd}\end{equation}
where the leftmost arrow is the symmetrization functor and the
rightmost diagonal arrow forgets the linear ordering and adds a basepoint.

We have the following localization result (see \autoref{sec:mainproof}):

\begin{Thm}\label{thm:localizationWF}
  The functor $\treerootsym\colon\symTrees\to \Gamma$
  exhibits $\Gamma$
  as an \inftycategorical{} localization of $\symTrees$
  at the set of collapse maps.
\end{Thm}

\begin{Rem}
  The functor $\treerootsym\colon\symTrees\to\finsetpop$ can be described as
  $\treerootsym\colon T\mapsto \lambda(T)\disjunion\set{\star}$,
  where $\lambda(T)$ is the set of leaves of a tree $T$.
  In this guise, it was introduced by \nameBM{}~\cite{BM17}.
\end{Rem}

\section{The localization theorem}

\label{sec:mainproof}
The following theorem expresses that the functor
$\treerootpl\colon\planarTrees\to\Delta$
(and its symmetric sibling $\treerootsym$)
is universal (in the \inftycategorical{} sense)
with the property of inverting the collapse maps in $\treerootpl$.

\begin{Thm}\label{thm:localizationW}
  For every \inftycat{} $\localC$,
  the functor
  $\treerootpl\colon \planarTrees\to\Delta$
  induces a fully faithful functor
  \begin{equation}
    \treerootpl^\star\colon \Fun(\Delta,\localC)
    \lra
    \Fun(\planarTrees,\localC)
  \end{equation}
  of \inftycats{} with essential
  image spanned by those functors $\planarTrees\to\localC$
  which map collapse maps $\corolla{}\to T$ to equivalences.
  The analogous statement holds for the functor
  $\treerootsym\colon \symTrees\to \Gamma$.
\end{Thm}

\begin{Cor}\label{cor:symtreecontractible}
  The categories
  $\planarTrees$ and $\symTrees$
  are weakly contractible
\end{Cor}
\begin{Prf}
  Clearly the categories $\Delta$ and $\Gamma$ are contractible
  because they have a terminal object and a zero object, respectively.
  Since the localization functors
  of \autoref{thm:localizationW} induce weak equivalences on classifying spaces,
  the result follows.
\end{Prf}

\begin{Rem}
  The weak contractibility of $\symTrees$
  (and implicitly of $\planarTrees$)
  was proved with a different method by
  \threenames{Ara}{Cisinski}{Moerdijk}~\cite{ACM2019}.
\end{Rem}

\subsection{The general situation}

Our strategy to prove \autoref{thm:localizationW} is to apply the following
general lemma which we will prove separately in
\autoref{subsection:proofoflocalization} below.
\begin{Lem}\label{lem:quillentheoremC}
  Let $\localL\colon W\to D$ be a functor of (ordinary) categories and for each
  $n\in D$ let $B_n\subset W_n$
  be a subcategory with the same objects of the weak fiber $W_n$ of $L$
  such that (with the notation of \autoref{rem:notationweakleftfibers} below)
  \begin{itemize}
  \item $B_n$ has an initial object $c_n$ and
  \item the inclusion $\nerve{B_n}\hra \overcat {\nerve W} n$ is cofinal.
  \end{itemize}
  Then for every \inftycat{} $\localC$, composition with $\localL$
  induces a fully faithful functor
  \begin{equation}\localL^\star\colon \Fun(\nerve {D}, \localC)\lra \Fun(\nerve {W},\localC)\end{equation}
  of \inftycats{} with the essential image spanned by those functors
  $\nerve{W}\to \localC$ which send all the edges of the form $c_n\to t$ in
  $\nerve{B_n}$ (for $n\in D$) to equivalences.
\end{Lem}
\begin{Rem}\label{rem:notationweakleftfibers}
  Recall that the weak fiber $W_n$ of
  $\localL\colon W\to D$
  is the category whose objects consist of an object $t\in W$ and an
  isomorphism $t\xra{\cong} n$ in $D$. The left fiber $\overcat W n\supset W_n$
  has objects $(t, f\colon t\to n)$ where $f$ is not required to be an
  isomorphism.
\end{Rem}

Let $\flexibleTrees{}$
be any one of the categories
$\planarTrees$ and $\symTrees$;
let $\treeroots$ be the corresponding functor
(among $\treerootpl$ and $\treerootsym$)
and denote its target
(which is either $\Delta$ or $\Gamma$)
by $\calD$.
For every object $\n\in\calD$ we denote by
$\overcat{\flexibleTrees{}}\n$ the left fiber, by $\flexibleTrees{}_\n$ the weak
fiber and by $\bp\n\subset\flexibleTrees{}_\n$ the subcategory of
$\flexibleTrees{}_\n$ with the same objects but only boundary preserving
morphisms. We shall now show that the functors $\treeroots$ satisfy the
requirements for \autoref{lem:quillentheoremC}, thus concluding the proof of
\autoref{thm:localizationW}.

\begin{Prop}\label{prop:situationC}
  Fix an object $\n\in\calD$.
  \begin{enumerate}
  \item The $n$\=/corolla $\corolla\n$ (together with any identification
    $\treeroots\corolla\n\xra{\cong}\n$) is an initial object in the category
    $\bp\n$.
  \item\label{inprop:bpreflective} The inclusion
    $\bp\n\subset\flexibleTrees{}_\n\hra \overcat{\flexibleTrees{}}\n$ has a
    left adjoint.\qedhere
  \end{enumerate}
\end{Prop}

\begin{dCor}
  The inclusion $\bp\n\hra\overcat{\flexibleTrees{}}\n$ is cofinal in the sense
  of Joyal~\cite[8.11]{Joyal2008}~\cite[Theorem 4.1.3.1]{Lurie2009}.
\end{dCor}

\begin{Prf}[of \autoref{prop:situationC}]
  The first statement is obvious.\\
  The functor $\overcat {\flexibleTrees{}} \n\to \bp\n$ is constructed as
  follows: Given an object $(T, f\colon\treeroots T=\m \to \n)$ we define the
  tree $T_f$ by glueing some corollas to $T$ along its outer edges (see also
  \autoref{figure:buildTftree}).
  We only describe this process explicitly for
  $\treeroots=\treerootpl$;
  the construction is analogous for $\treerootsym$.
  \begin{itemize}
  \item To a leaf of $T$ corresponding to the minimal edge $\{j-1, j\}\hra \m$
    we glue a corolla $\corolla[f]{j-1,j}$ (of arity $f(j)-f(j-1)$) with leaves
    $\{i-1, i\}$ for $f(j-1)<i\leq f(j)$ (this might be a $0$\=/corolla if
    $f(j-1)=f(j)$).
  \item To the root (corresponding to the maximal edge $\{0, m\}\hra \m$) we
    glue a corolla $\corolla[f]\max$ with leaves
    \begin{equation}\{0,1\}, \{1,2\}, \dots , \{f(0)-1, f(0)\}, \{f(0), f(m)\},
      \{f(m),f(m)+1\}, \dots, \{n-1,n\}\end{equation} along the special leaf $\{f(0),
    f(m)\}$ of $\corolla[f]\max$.
  \end{itemize}

\begin{figure}[t]
  \centering {%
    \setlength{\fboxsep}{10pt}
    \fbox{\includegraphics[scale=0.8]{\figuresWhere/buildtree_Tf}} }%
  \caption{The construction of the tree $T_f$ in the case
    $\treeroots=\treerootpl$. The little arrows decorate the roots of the
    various trees. Forgetting the root and/or the plane embedding describes the
    analogous construction in the cases $\treeroots=\treerootcyc, \treerootsym,
    \treerootabstract$}
  \label{figure:buildTftree}

\end{figure}

The adjunction unit at $(T,f)$ is the inclusion $T\hra T_f$ which we denote by
$f_T$. We need to prove that given a morphism of trees $\alpha\colon T\to S$
over $f\colon \m\to\n$ there is a unique factorization
$T\xra{f_T}T_f\xra{\morbp\alpha}S$ with $\morbp\alpha$ in $\bp\n$. We have no
other choice than to define $\morbp\alpha$ as $\alpha$ on the subtree $T\hra
T_f$ and to make it the identity on the boundary; hence uniqueness is clear. It
is straightforward to verify that this map of trees is indeed well defined.
\end{Prf}

\subsection{Proof of the key lemma}

\label{subsection:proofoflocalization}

This section is devoted to the proof of \autoref{lem:quillentheoremC} 

Let $M$ be defined as the Grothendieck construction of the functor $\simplex
1\to \Cat$ which parameterizes the functor $L\colon W\to D$. Explicitly, an
object in $M$ is either an object $t\in W$ or an object $n\in D$; for $s,t\in W$
and $m,n\in D$ we put $M(t,s)=W(t,s)$ and $M(n,m)=D(n,m)$ and $M(t,n)=D(Lt, n)$
and $M(n,t)=\emptyset$. We have a factorization $\localL\colon W\hra M
\xra{\localLp}D$ where the first arrow is the obvious fully faithful inclusion
and the second arrow has a fully faithful right adjoint $D\hra M$. We identify
$D$ with its image in $M$ and we denote by $\adjunit\colon \Id_M\to \localLp$
the unit of the adjunction $\localLp\colon M\ladjarrows D$; it is an isomorphism
(in fact the identity) at exactly those objects in $M$ that belong to
$D$.\footnote{
  The components $\adjunit_t\colon t\to \localLp t$ of the
  adjunction are precisely the coCartesian morphisms of the coCartesian
  fibration $M\to \simplex 1$.
}
We deal with the two components of $\localL\colon W\hra M\ladjarrows D$
individually by using standard techniques from Higher Topos
Theory~\cite{Lurie2009}. \autoref{lem:quillentheoremC} is a direct
consequence of \autoref{lem:secondhalfthmC} and \autoref{lem:firsthalfthmC}
below.

\begin{Rem}
  For each $n\in D$ the forgetful functor $B_n\subset W_n \to W$ extends to a
  functor $\rcone{B_n}\hra M$ by sending the new vertex $v$ to $n$ and the new
  arrow $(t,f)\to v$ (for $(t,f)\in B_n$) to the arrow $f\colon t\to n$ of $M$.
\end{Rem}

Fix an \inftycat{} $\localC$. We recall the following result.

\begin{Lem}~\cite[Proposition~5.2.7.12]{Lurie2009}
  \label{prop:HTT5.2.7.12}
  Let $\localLp \colon \calM\to \calD$ be a reflective localization functor of
  \inftycats{} (i.e.\ $\localLp $ has a fully faithful right adjoint)
  and let $\localC$ be another \inftycat{}. Then composition with
  $\localLp $ induces a fully faithful functor
  \begin{equation}\Fun(\calD,\localC)\lra \Fun(\calM,\localC)\end{equation}
  with essential image consisting of those functors that map an edge $f$ in
  $\calM$ to an equivalence in $\localC$ provided that $\localLp f$ is an
  equivalence in $\calD$.
\end{Lem}

\begin{Lem}\label{lem:descriptionsKplus}
  Let $F\colon \nerve M \to \localC$ be a functor of \inftycats{}. The
  following are equivalent:
  \begin{enumerate}
  \item\label{lem:item:alllocalizingmaps} For every edge $f$ in $\nerve M$, if
    $\localLp f$ is an equivalence in $D$ then $F f$ is an equivalence in
    $\localC$.
  \item\label{lem:item:thosemapsinBncone} For every $n\in D$, the functor $F$
    maps all edges in $\rcone{\nerve{B_n}}$ to equivalences in $\localC$.
  \item\label{lem:item:onlyadjunctionmaps} $F$ sends every component
    $\adjunit_t\colon t\to \localLp t$ of the unit to an equivalence in
    $\localC$.
  \end{enumerate}
  We denote by $K^+$ the full subcategory of $\Fun(\nerve{M},\localC)$ spanned
  by such functors.
\end{Lem}

\begin{Prf}
  Clearly \ref{lem:item:alllocalizingmaps} implies
  \ref{lem:item:thosemapsinBncone} because $\localLp f$ is an isomorphism for
  each edge $f$ of $\rcone{\nerve{B_n}}$.
  Condition \ref{lem:item:thosemapsinBncone}
  implies \ref{lem:item:onlyadjunctionmaps} because for each $t\in W$,
  the edge $\adjunit_t$ appears in the cone $\rcone{B_n}\to M$
  as the structure map
  over $(t,\localL t\xrightarrow{=}\localL t\eqqcolon n)\in B_n$
  (and for $t\in D$ the component $\eta_t$ is the identity;
  hence that case is automatic).\\
  Observe that if $f\colon t\to s$ is a morphism in $M$ then we have a
  commutative naturality square
  \begin{equation}\begin{tikzcd}
      t\ar[d, "f"]\ar[r,"\adjunit_t"]& \localLp t\ar[d,"\localLp f"]\\
      s\ar[r,"\adjunit_s"]& \localLp s
    \end{tikzcd}\end{equation} Hence \ref{lem:item:onlyadjunctionmaps} implies
  \ref{lem:item:alllocalizingmaps} by the two-out-of-three property for
  equivalences in $\localC$.
\end{Prf}

\begin{dCor}\label{lem:secondhalfthmC}
  Composition with the functor $\localLp \colon M\to D$ induces a fully faithful
  functor $\Fun(\nerve{D},\localC)\hra\Fun(\nerve{M},\localC)$ with essential
  image $K^+$.
\end{dCor}

Let us recall the following result.

\begin{Lem}~\cite[Proposition~4.3.1.12]{Lurie2009}
  \label{prop:HTT4.3.1.12}
  Let $\localC$ be an \inftycat{} and let $\ol F\colon \rcone
  B\to\localC$ be a diagram where $B$ is a weakly contractible simplicial set
  and $\ol F$ carries each edge of $B$ to an equivalence in $\localC$. Then $\ol
  F$ is a colimit diagram in $\localC$ if and only if it carries every edge in
  $\rcone B$ to an equivalence in $\localC$.
\end{Lem}

\begin{Lem}\label{lem:descriptionsK}
  Let $F\colon \nerve W \to \localC$ be a functor. The following are equivalent:
  \begin{enumerate}
  \item\label{lem:item:locallyconstantkanextension} The functor $F$ admits a
    left Kan extension along $W\hra M$ and the resulting functor $\nerve M \to
    \localC$ lies in $K^+$.
  \item\label{lem:item:allbpmaps} For every $n\in D$ the functor $F$ maps every
    edge of $\nerve{B_n}$ to an equivalence in $\localC$.
  \item\label{lem:item:onlycollapsemaps} For every $n\in D$ and every $t\in B_n$
    the functor $F$ maps the unique edge $c_n\to t$ in $\nerve{B_n}$ to an
    equivalence in $\localC$.
  \end{enumerate}
  We denote by $K$ the full subcategory of $\Fun(\nerve{W},\localC)$ spanned by
  such functors.
\end{Lem}

\begin{Prf}
  The equivalence between \ref{lem:item:allbpmaps} and
  \ref{lem:item:onlycollapsemaps} is obvious because $c_n$ is an initial element
  in $B_n$. Using description~\ref{lem:item:thosemapsinBncone} of
  \autoref{lem:descriptionsKplus} it is clear that
  \ref{lem:item:locallyconstantkanextension} implies \ref{lem:item:allbpmaps}.

  Let us prove the converse: By the pointwise construction of Kan
  extensions~\cite[Lemma 4.3.2.13]{Lurie2009}, a left Kan extension
  of $F$ along $W\hra M$ can be assembled from colimit cones for the diagrams
  $\overcat {\nerve W} n\to \nerve W\xra{F}\localC$ (for $n\in D$). Recall that
  $B_n\hra \overcat W n$ is cofinal, hence we can reduce to finding colimits for
  the diagrams $\nerve {B_n}\hra \nerve{\overcat W n}\to
  \nerve{W}\xra{F}\localC$. All edges of these diagrams are equivalences by
  condition~\ref{lem:item:thosemapsinBncone} and $\nerve {B_n}$ is contractible
  (because $B_n$ has an initial element). Therefore by
  \autoref{prop:HTT4.3.1.12} these colimits exists and the corresponding colimit
  cones $\rcone{\nerve{B_n}}\to\localC$ map all edges to equivalences in
  $\localC$, thus verifying condition~\ref{lem:item:thosemapsinBncone} of
  \autoref{lem:descriptionsKplus}.
\end{Prf}

\newcommand\localH{H} Fix the following notation:
\begin{itemize}
\item Denote by $\localH^+$ the full subcategory of $\Fun(\nerve{M},\localC)$
  spanned by those functors which are the left Kan extension of their
  restriction to $W\subset M$.
\item Denote by $\localH$ the full subcategory of $\Fun(\nerve{W}, \localC)$
  spanned by those functors which admit a left Kan extension along $W\hra M$.
\end{itemize}

Recall the following result.

\begin{Lem}~\cite[Proposition~4.3.2.15]{Lurie2009}
  The restriction functor along $\nerve W\hra\nerve M$ is a trivial fibration
  $\localH^+\to \localH$ of simplicial sets.
\end{Lem}

\begin{Lem}\label{lem:pullbackconstantcolimits}
  We have inclusions $K^+\subset \localH^+$ and $K\subset \localH$ and a
  pullback square
  \begin{equation}\begin{tikzcd}
      K^+\ar[d]\ar[r, hookrightarrow]& \localH^+\ar[d]\\
      K\ar[r, hookrightarrow]&\localH
    \end{tikzcd}\end{equation} of simplicial sets with vertical arrows given by restriction
  along $W\hra M$.
\end{Lem}

\begin{Prf}
  This follows directly from \autoref{lem:descriptionsKplus} and
  \autoref{lem:descriptionsK} \end{Prf}

Since trivial fibrations of simplicial sets are stable under pullbacks we
obtain:

\begin{dCor}\label{lem:firsthalfthmC}
  The restriction functor along the inclusion $W\hra M$ is a trivial fibration
  $K^+\to K$ of simplicial sets.
\end{dCor}

This concludes the proof of \autoref{lem:quillentheoremC} and therefore of
\autoref{thm:localizationW}.

\section{Applications}

Consider the category
$\sSet\coloneqq [\Dop, \Set]$
of simplicial sets equipped
with the Kan--Quillen
left proper combinatorial simplicial model structure~\cite{Quillen1967}.
Denote by
$\Spaces\coloneqq \snerve{\sSet^\fibcofibdec}$
the corresponding \introduce{\inftycat{} of spaces}
obtained as the simplicial nerve of the subcategory
of fibrant-cofibrant objects~\cite[Definition~1.2.16.1]{Lurie2009}.
A dendroidal (\resp simplicial) object in $\Spaces$
is called a dendroidal (\resp simplicial) space.

\subsection{\twoSegal{} simplicial objects and Segal dendroidal objects}

\label{sec:SegalvsSegal}

In this section we compare the dendroidal Segal conditions due to
\nameCM{}~\cite{CM13}
and the simplicial \twoSegal{} conditions due to Dyckerhoff and
Kapranov~\cite{DyckerhoffKapranov2012}.

\begin{Def}\label{def:SegalCM}~\cite[Definition 2.2]{CM13}
  The
  \introduce{Segal core} of a tree $\degtree\neq T\in\symTrees$ is the union
  \begin{equation}\Segalcore T\coloneqq \bigcup\limits_{v}\symTrees[\corolla {n(v)}]\end{equation}
  where $v$ runs over all vertices of $T$ and $\corolla{n(v)}\hra T$ denotes the
  subtree with vertex $v$. We use the convention $\Segalcore\degtree \coloneqq
  \symTrees[\degtree]$ for the trivial tree.

  A symmetric dendroidal space $\localX\colon \symTreesop\to\Spaces$ is
  \introduce{Segal} if for any tree $T\in \symTrees$ the map
  \begin{equation}\localX_T=\Hom(\symTrees(T),\localX)\lra \Hom(\Segalcore T, \localX)\end{equation}
  is a trivial fibration.
\end{Def}

We adapt this definition as follows.

\begin{Def}\label{def:Segalgrafting}
  A dendroidal object $\localX\colon\planarTreesop\to \localC$ in some
  \inftycat{} $\localC$ is called \introduce{Segal} if $\localX$ sends
  the diagram
  \begin{equation}\label{equ:indeftreesegal}\begin{tikzcd}
      T\ar[r,hookleftarrow]&T_2\\
      T_1\ar[u,hookrightarrow]\ar[r,hookleftarrow]&e\ar[u,hookrightarrow]
    \end{tikzcd}\end{equation}
  to a pullback square in $\localC$ whenever the tree $T\in\planarTrees$ arises
  by grafting two trees $T_1$ and $T_2$ along a common edge $e$.
\end{Def}

\begin{Rem}
  Clearly \autoref{def:SegalCM} and \autoref{def:Segalgrafting} make sense,
  mutatis mutandis, for symmetric dendroidal objects.
  Another way of saying this is that a symmetric dendroidal object is Segal
  if and only if the underlying dendroidal object is Segal.
\end{Rem}

\begin{Rem}
  If a tree $T$ arises by grafting two trees $T_1$ and $T_2$ along a common edge
  $e$ then clearly $\Segalcore T = \Segalcore {T_1}\pushout{e}\Segalcore{T_2}$.
  By successively decomposing a tree along its inner edges we therefore see that
  \autoref{def:SegalCM} and \autoref{def:Segalgrafting} agree for dendroidal
  objects in the \inftycat{} $\Spaces$ of spaces.
\end{Rem}

The importance of the dendroidal Segal conditions
is highlighted by the following result,
which has an obvious analog for \non{symmetric} operads and dendroidal sets.

\begin{Prop}\label{prop:CMSegaloperads}~\cite[Corollary 2.6]{CM13}
  The symmetric dendroidal nerve functor
  \begin{equation}\dnerve{}\colon \symOperads\lra\dSet\end{equation}
  is fully faithful and the essential image consists precisely of the Segal
  symmetric dendroidal sets.
\end{Prop}

\begin{Def}\label{def:unital2Segal}~\cite[Proposition 2.3.2]{DyckerhoffKapranov2012}
  A simplicial object
  $\localX\colon\Dop\to\localC$
  in some \inftycat{} $\localC$
  is called \introduce{\twoSegal{}}
  (or \introduce{unital \twoSegal{}})
  if for each $0\leq i\leq j\leq m$
  it maps the square
  \begin{equation}\label{equ:indef2segal}\begin{tikzcd}
      \{0,\dots, m\}\ar[r,hookleftarrow]& \{i,\dots, j\}\\
      \{0,\dots,i,j, \dots m\}\ar[u]\ar[r,hookleftarrow]&\{i,j\}\ar[u]
    \end{tikzcd}\end{equation}
  in $\Delta$ to a pullback square square in $\localC$.
\end{Def}
\begin{Rem}
  We always interpret the elements $i$ and $j$ in the lower row of
  Diagram~\ref{equ:indef2segal} as distinct; thus in the case $i=j$ the vertical
  arrows are codegeneracy maps.
\end{Rem}

\begin{Rem}
  The original definition of \twoSegal{} objects
  only includes the case $i\neq j$ of \eqref{equ:indef2segal};
  the condition for $i=j$ is called \buzzword{unitality}.
  Since unitality is now known to be redundant~\cite{FGKUPW2019},
  we drop that adjective entirely.
\end{Rem}

\begin{Lem}\label{lem:comparisonofSegals}
  A simplicial object $\localX\colon\Dop\to \localC$ in some
  \inftycat{} $\localC$ is \twoSegal{} if and only if the composition
  $\treerootpl^\star
  \localX\colon\planarTreesop\xra{\treerootpl}\Dop\xra{\localX}\localC$
  is a Segal dendroidal object.
\end{Lem}

\begin{Prf}
  Let $T=T_1\graft[e] T_2$ be a grafting of trees where $e$ is the root of $T_2$
  and a leaf of $T_1$. Put $\m\coloneqq \treerootpl T$. Applying $\treerootpl$
  to the inclusion $e\hra T$ defines a map ${\numD 1}=\treerootpl e\xra{f} \m$, so we
  can define $i\coloneqq f(0)$ and $j\coloneqq f(1)$. It is easy to see that
  with this notation $\treerootpl$ sends Diagram~\eqref{equ:indeftreesegal} to
  Diagram~\eqref{equ:indef2segal} and that every instance of
  Diagram~\eqref{equ:indef2segal} arises this way.
\end{Prf}

\subsection{Segal simplicial objects and covariantly fibrant dendroidal objects}

\label{sec:Segalvscovariant}
Recall that a simplicial object $\localX\colon \Dop\to \localC$ in some
\inftycat{} $\localC$ is called \introduce{reduced Segal} if
$\localX_\n\xra{\simeq} \localX_{\numD 1}^n$ via the inert maps
$\{i-1,i\}\hra \numD{n}$ in $\Delta$
(in particular $\localX_{\numD 0}$ is a terminal object in $\localC$).
A similar condition makes sense when replacing $\Delta$ by
$\Gamma\coloneqq\finsetpop$; such functors $\localX\colon \Gamma^\op\to\localC$
were introduced (in the case $\localC\coloneqq \Spaces$) by
Segal~\cite{Segal1974} under the name
\buzzword{special $\Gamma$-spaces}.

\begin{Def}~\cite{BM17}
  \label{def:covariantly-fibrant}
  A dendroidal object
  $\localX\colon\planarTreesop\to\localC$
  (or $\localX\colon\symTreesop\to\localC$)
  is \introduce{covariantly fibrant}
  if for each $n$-ary tree $T$
  the inclusion of its leaves $l_1,\dots, l_n$,
  induces an equivalence
  $\localX_T\xra{\simeq} \prod_{i=1}^n \localX_{l_i}$.
\end{Def}

It is clear from the definitions that
\begin{itemize}
\item a simplicial object $\localX$ in $\localC$ is reduced Segal if and only if
  $\treerootpl^\star\localX$ is covariantly fibrant,
\item every covariantly fibrant $\localX\colon\planarTreesop\to \localC$ maps
  collapse maps to equivalences.
\end{itemize}
(And similarly for the symmetric case.)
Therefore \autoref{thm:localizationW}
immediately implies the following result, proved by
\nameBM{}~\cite[Theorem 1.1]{BM17}
for $\localC=\Spaces$ in the language of model categories.

\begin{dCor}\label{cor:Segalvscovariantfibration}
  For every \inftycat{} $\localC$, the functor $\treerootpl$ (\resp
  $\treerootsym$) induces an equivalence of \inftycats{} between
  \begin{itemize}
  \item reduced Segal simplicial (\resp $\Gamma$-) objects in $\localC$
  \item covariantly fibrant plane (\resp symmetric) dendroidal objects in
    $\localC$.\qedhere
  \end{itemize}
\end{dCor}

\subsection{\twoSegal{} simplicial sets and invertible operads}

\label{sec:invertibleoperads}

\begin{Def}~\cite[Definition 3.6.7]{DyckerhoffKapranov2012}
  \label{def:invertible-operad}
  An operad $\localO$ is called \introduce{invertible}
  if the unit map
  \eqref{eq:operad-big-unit}
  and all the composition maps
  \eqref{eq:operad-big-composition}
  are invertible.
\end{Def}

\begin{Rem}
  It follows from
  \autoref{rem:circle-i-compositions}
  that an operad is invertible
  if and only if
  the unit map
  \eqref{eq:operad-big-unit}
  and all $\circleis[j]{i}$-compositions
  \eqref{eq:circle-i-compositions}
  are invertible.
\end{Rem}

\begin{Rem}
  It follows from the condition on the unit map
  that if an operad is invertible then its underlying category is discrete,
  \ie has only identity arrows.
\end{Rem}

\begin{Prop}~\cite[Theorem 3.6.8]{DyckerhoffKapranov2012}\label{prop:thmDKoperads}
  Fix a set $B$ of colors. Then there is an equivalence of categories between
  invertible $B$-colored operads and \twoSegal{} simplicial sets
  $\localX\colon\Dop\to\Set$ with $\localX_{\numD 1}=B$.
\end{Prop}

We can characterize invertibility of an operad in terms of its dendroidal nerve.
\begin{Lem}\label{lem:charinvertibleoperads}
  Let $\localO$ be an operad and let $\dnerve\localO\colon
  \planarTreesop\to\Set$ be its dendroidal nerve. The following are equivalent:
  \begin{enumerate}
  \item\label{inprop:nerveinvertsbp} The dendroidal set $\dnerve\localO$ maps
    all boundary preserving maps to isomorphisms.
  \item\label{inprop:nerveinvertscollapse} The dendroidal set $\dnerve\localO$
    is invertible, i.e.\ it inverts all collapse maps.
  \item\label{inprop:Oisinvertible} The operad $\localO$ is invertible.
  \qedhere
  \end{enumerate}
\end{Lem}

\begin{Prf}
  If $\alpha\colon T\to S$ is boundary preserving, then clearly the collapse map
  for $S$ factors through the collapse map for $T$ as $\corolla {}\to
  T\xra{\alpha}S$. Hence \ref{inprop:nerveinvertsbp} and
  \ref{inprop:nerveinvertscollapse} are equivalent by the
  $2$-out-of-$3$-property for isomorphisms.

  The unit map~\eqref{eq:operad-big-unit}
  in \autoref{def:coloredoperad}
  is precisely the image under $\dnerve \localO$ of the collapse map
  $\corolla 1\to \eta$.
  Taking the coproduct over all the composition maps for
  fixed $k, n_1,\dots, n_k\in\BN$ yields
  (putting $n \coloneqq \sum_{i=1}^k n_i$)
  precisely the image of the collapse map
  $\corolla n \to T_k^{n_1,\dots,n_k}$,
  where $T_k^{n_1,\dots,n_k}$ is tree obtained by glueing
  (for all $0\leq i\leq k$) the corolla $\corolla {n_i}$ to the $i$-th leaf of
  the corolla $\corolla k$.
  Hence \ref{inprop:nerveinvertscollapse} implies \ref{inprop:Oisinvertible}.
  The converse holds because every
  \roughly{generalized composition map} represented by a collapse map
  $\corolla{}\to T$ can be written as the composition of unit and composition
  maps as in \autoref{def:coloredoperad}.
\end{Prf}

Using
\begin{itemize}
\item the characterization of operads as Segal dendroidal sets (the
  \non{symmetric} analogue of \autoref{prop:CMSegaloperads}),
\item the characterization of invertible operads
  (\autoref{lem:charinvertibleoperads}),
\item our main result (\autoref{thm:localizationWD}) in the case $\localC=\Set$
  and
\item the corresponcence between Segal dendroidal objects and \twoSegal{}
  simplicial objects (\autoref{lem:comparisonofSegals})
\end{itemize}
we recover the following more elegant version of \autoref{prop:thmDKoperads}.

\begin{dCor}
  The composition $\sSet\xra{\treerootpl^\star}\pldSet\lra \Operads$ restricts
  to an equivalence of categories between the full subcategories of \twoSegal{}
  simplicial sets on one side and invertible operads on the other.
\end{dCor}

\newcommand\zerobj{0}
\newcommand\exabcat{A}
\newcommand\Scstr{\rmS}
\newcommand\ScstrofA{\Scstr(\exabcat)}
\newcommand\ScAop[1]{\ScstrofA_{#1}}
\newcommand\localx{x}
\newcommand\nn{n}
\newcommand\kk{k}
\newcommand\xux[2]{{_{#1}\localx_{#2}}}
\newcommand\locala{i}
\newcommand\localb{j}

Before moving on,
we discuss some examples of invertible operads.

\begin{Expl}[Waldhausen's \Sconstruction{}~\cite{Waldhausen1985}]
  \label{expl:WaldhausenS}
  Let $\exabcat$ be an abelian category\footnote{
    Waldhausen's \Sconstruction{} applies in much greater generality;
    we restrict to abelian categories for simplicity.}.
  Consider the following operad $\ScstrofA$:
  \begin{itemize}
  \item
    The colors of $\ScstrofA$ are
    the objects of $\exabcat$ (up to isomorphism).
  \item
    The \ary{2} operations of $\ScstrofA$ are
    short exact sequences
    \begin{equation}
      \cdsquareNA[bC]{\xux01}{\xux12}{\zerobj}{\xux02}
    \end{equation}
    (up to isomorphism)
    each of which is viewed as a \ary{2} operation $(\xux01,\xux12)\to\xux02$.
  \item
    More generally, the \ary{\nn} operations
    $
    (\xux01,\xux12,\dots,\xux{\nn-1}\nn)\lra\xux0\nn
    $
    of $\ScstrofA$ are diagrams
    \begin{equation}
      \label{eq:flag-operation}
      \newcommand\squu{\ar[r]\ar[d]\isbiCartesian}
      \newcommand\sqss{\ar[r]\ar[d]}
      \begin{tikzcd}
        \xux01\squu
        &\xux02\sqss
        &\cdots\ar[r]
        &\xux0{\nn-2}\squu
        &\xux0{\nn-1}\squu
        &\xux0\nn\ar[d]\\
        \zerobj\ar[r]
        &\xux12 \sqss
        &\cdots\ar[r]
        &\xux0{\nn-2}\squu
        &\xux0{\nn-1}\squu
        &\xux0\nn\ar[d]\\
        &\zerobj\ar[r] 
        &\cdots\ar[r]
        &\xux0{\nn-2}\sqss
        &\xux0{\nn-1}\sqss
        &\xux0\nn\ar[d]\\
        &&\ddots
        &\vdots\ar[d]
        &\vdots\ar[d]
        &\vdots\ar[d]\\
        &&&\zerobj\ar[r]
        &\xux{\nn-2}{\nn-1}\squu
        &\xux{\nn-2}\nn\ar[d]\\
        &&&&\zerobj\ar[r]
        &\xux{\nn-1}\nn\\
      \end{tikzcd}
    \end{equation}
    in $\exabcat$
    (up to isomorphism),
    where each square is required to be biCartesian,
    \ie both a pushout and a pullback.
  \item
    The $\circleis[\localb]{\locala}$ composition of an operation
    \begin{equation}
      f\colon
      (\xux\locala{\locala+1},\dots,\xux{\localb-1}\localb)
      \lra\xux{\locala}\localb
    \end{equation}
    with an operation
    \begin{equation}
      g\colon
      (\xux01,\dots,\xux\locala\localb,\dots,\xux{\nn-1}\nn)
      \lra\xux0\nn
    \end{equation}
    (for $0\leq\locala\leq\localb\leq\nn$)
    is the operation
    \begin{equation}
      (g\ocomp[\locala+1]f)\colon
      (\xux01,\xux12,\dots,\xux{\nn-1}\nn)
      \lra\xux0\nn
    \end{equation}
    whose associated diagram \eqref{eq:flag-operation}
    is uniquely characterized
    by the fact that it
    extends the corresponding diagrams for $f$ and $g$.
  \end{itemize}
  It is not hard to verify that $\ScstrofA$ is a well defined operad;
  it is invertible because,
  for each fixed $0\leq\locala\leq\localb\leq\nn$,
  each operation \eqref{eq:flag-operation}
  arises as the composition
  $g\ocomp[\locala+1]f$
  for a unique pair of operations $(f,g)$ as above.
  Under the equivalence of \autoref{thm:intro-Gamma-Lambda}
  this operad corresponds to Waldhausen's \Sconstruction{}
  which is the \twoSegal{} simplicial set
  $\ScstrofA\colon\Dop\to\Set$
  that maps $\numD\nn\in\Delta$
  to the set of
  isomorphism classes of
  diagrams \eqref{eq:flag-operation}
  with \twoptions{face}{degneracy} maps given
  by simultaneously \twoptions{omitting}{duplicating} rows and columns.
  If instead of working up to isomorphism
  we keep track of those isomorphisms,
  we get an \twoptions{invertible operad}{\twoSegal{} object} in groupoids
  rather than sets.
\end{Expl}

\begin{Rem}
  Let $\localX$ be an invertible Segal dendroidal object.
  Let $T$ be the closed $\nn$\=/corolla
  (i.e. the grafting of $\nn$ many $0$\=/corollas on top of a $\nn$\=/corolla).
  We have two maps
  \begin{equation}\localX(\corolla{0})\xla{\simeq}
    \localX(T)\xra{\simeq}\localX(\corolla{n})\fiberproduct{\localX(\degtree)^n}\localX(\corolla{0})^n
  \end{equation}
  which are equivalences by invertibility and the Segal conditions respectively.
  In the example where $\localX=\ScstrofA$ is the Waldhausen \Sconstruction{}
  of an abelian category $\exabcat$,
  the groupoid $\localX(\corolla{0})\simeq \set{\zerobj}$ is trivial,
  hence this condition says precisely that
  a flag \eqref{eq:flag-operation} of length $\nn$
  with trivial \sub{quotients} is trivial.
  Note, however,
  that in general a flag is not determined by its \sub{quotients},
  which would be the Segal condition
  $\localX(\corolla{n})\xra{\simeq}\localX(\degtree)^n$.
\end{Rem}

\begin{Expl}[\Eoperad{\kk}s]
  The commutative operad $\Ealg\infty$
  (viewed as a symmetric operad)
  has a contractible space of operations in each degree,
  hence is invertible for trivial reasons;
  it corresponds to the constant $\Gamma$-space on a point.
  Its underlying \non{symmetric} operad
  is the associative operad
  which is invertible and
  corresponds to the constant simplicial space on a point.
  For all other $1\leq\kk<\infty$,
  the operad $\Ealg\kk$ of little $\kk$-cubes
  is easily seen to \emph{not} be invertible.
\end{Expl}

\begin{Expl}
  \label{expl:monoid-as-operad}
  \newcommand\exmonoid{M}
  \newcommand\exelofmon{m}
  \newcommand\mum[2]{{_{#1}\exelofmon_{#2}}}
  \newcommand\locali{i}
  Each monoid $\exmonoid$
  (multiplicatively written)
  gives rise to an invertible operad $\nerve{\exmonoid}$ as follows:
  The set of colors is $\exmonoid$.
  The set of \ary{\nn} operations is $\exmonoid^\nn$,
  where each tuple
  $(\mum01,\dots,\mum{\nn-1}{\nn})\in\exmonoid^\nn$
  is viewed as an operation
  \begin{equation}
    (\mum01,\dots,\mum{\nn-1}{\nn})
    \lra
    \mum01\cdots\mum{\nn-1}{\nn}
    \eqqcolon
    \mum0{\nn}
  \end{equation}
  and is, for each $0\leq\locala\leq\localb\leq\nn$,
  the $\circleis[\localb]{\locala}$-composition
  of
  \begin{equation}
    (\mum{\locala}{\locala+1},\dots,\mum{\localb-1}{\localb})
    \lra
    \mum{\locala}{\localb}
  \end{equation}
  and
  \begin{equation}
    (
    \mum01,
    \dots,
    \mum{\locala-1}{\locala},
    \mum{\locala}{\localb},
    \mum{\localb}{\localb+1},
    \dots,
    \mum{\nn-1}{\nn}
    )
    \lra
    \mum0{\nn}.
  \end{equation}
  If $\exmonoid$ is abelian
  then the operad $\nerve{\exmonoid}$
  can be canonically enhanced to a symmetric operad.
  Under the equivalence of \autoref{thm:intro-Gamma-Lambda},
  the operad $\nerve{\exmonoid}$
  corresponds to the nerve
  $\nerve{\exmonoid}\colon\Dop\to\Set$
  which is not just \twoSegal{} but Segal.
  
  This example can be categorified
  to interpret each monoidal \inftygrpd{}
  as an invertible \inftyoperad{};
  see \autoref{expl:monoidal-grpd-cov-fibrant}
  and \autoref{rem:operads-monoidal-grpds}.
\end{Expl}

\subsection{\twoSegal{} simplicial spaces and invertible \inftyoperad{}s}

\label{sec:invertibleinfoperads}

As a direct consequence of \autoref{thm:localizationWD} and
\autoref{lem:comparisonofSegals} we obtain the following comparison result.

\begin{dCor}\label{cor:Segaldendrsimplspaces}
  Composition with $\treerootpl\colon\planarTrees\to\Delta$ induces an
  equivalence between the \inftycat{} of \twoSegal{} simplicial spaces
  and the \inftycat{} of invertible Segal dendroidal spaces.
\end{dCor}

The goal of this \autoref{sec:invertibleinfoperads} is to give an
interpretation of this result by identifying the \inftycat{} of
invertible Segal dendroidal spaces as a full subcategory of the
\inftycat{} of complete Segal dendroidal spaces. We treat the latter as
a model for (\non{symmetric}) \inftyoperad{}s (in analogy to results due to
\nameCM{}~\cite{CM13}
in the symmetric case) so that we can rephrase
\autoref{cor:invinfopsSegalspaces} as follows:

\begin{dCor}\label{cor:invinfopsSegalspaces}
  Composition with $\treerootpl\colon\planarTrees\to\Delta$ induces an
  equivalence between the \inftycat{} of \twoSegal{} simplicial spaces
  and the \inftycat{} of invertible (\non{symmetric}) \inftyoperad{}s.
\end{dCor}

\begin{Expl}
  \label{expl:monoidal-grpd-cov-fibrant}
  \newcommand{\exmongrpd}{\calM}
  \newcommand{\exmonoperad}{\localO_\calM}
  Every monoidal category $(\exmongrpd,\otimes)$
  gives rise to an operad $\exmonoperad$
  in groupoids:
  Its groupoid of colors
  $\exmonoperad(\degtree)\coloneqq\grpdcore{\exmongrpd}$
  is the groupoid core
  of $\exmongrpd$
  and its groupoid of \ary{1} operations
  is the groupoid
  $
  \exmonoperad(1)\coloneqq
  \grpdcore{\Fun(\simplex{1},\exmongrpd)}
  $
  of arrows in $\exmongrpd$.
  The groupoid $\exmonoperad(n)$ of \ary{n} operations
  is the groupoid of arrows
  $\bullet_1\otimes\dots\otimes\bullet_n\to\bullet$,
  \ie the pullback
  \begin{equation}
    \label{eq:pullback-for-tensor-arrows}
    \cdsquare[C]
    {\exmonoperad(n)}
    {\exmonoperad(1)}
    {\exmonoperad(\eta)^n}
    {\exmonoperad(\eta)}
    {}
    {}
    {\mathrm{s}}
    {\otimes}
  \end{equation}
  Composition in the operad $\exmonoperad$
  is induced by composition of arrows in $\exmongrpd$.
  The operad $\exmonoperad$ is invertible
  if and only if
  all arrows in the underlying category $\exmongrpd$ are invertible,
  \ie if and only if
  $\exmongrpd$ is a monoidal \emph{groupoid}.
  In this case,
  the right vertical map in \eqref{eq:pullback-for-tensor-arrows}---%
  which sends each arrow to its source---%
  is an equivalence;
  hence the same is true for the left vertical map.
  This amounts to saying that,
  viewed as a Segal dendroidal groupoid,
  $\exmonoperad$ is covariantly fibrant.

  Under the equivalence of \autoref{cor:Segalvscovariantfibration},
  the operad $\exmonoperad$
  corresponds to the complete Segal simplicial space
  obtained by interpreting $\exmongrpd$
  as an \inftycat{} with a single object,
  $\exmongrpd$ as its space of arrows and composition given by $\otimes$.
  This generalizes
  \autoref{expl:monoid-as-operad},
  where the monoidal groupoid $\exmongrpd$ is discrete.
\end{Expl}

\begin{Rem}
  \label{rem:operads-monoidal-grpds}
  In view of \autoref{expl:monoidal-grpd-cov-fibrant}
  and considering that complete reduced Segal simplicial spaces
  are a model for monoidal \inftygrpd{}s\footnote{
    For instance,
    Lurie~\cite[Definition~4.1.3.6]{Lurie2017} defines
    (\non{symmetric}) monoidal \inftycats{}
    as those coCartesian fibrations over $\Dop$
    which under the \twoptions{straightening}{unstraightening} equivalence
    correspond to reduced Segal simplicial \inftycats{};
    monoidal \inftygrpd{}s are then precisely
    those that take values in \inftygrpd{}s rather than \inftycats{}.
  },
  \autoref{cor:Segalvscovariantfibration}
  allows us to interpret
  \qquote{being covariantly fibrant}
  as the property which characterizes those \inftyoperads{}
  which are monoidal \inftygrpd{}s.
\end{Rem}

The theory of complete Segal dendroidal spaces was developed by
\nameCM{}~\cite{CM13}
and spelled out in detail for \emph{symmetric} dendroidal spaces. They prove
that complete Segal symmetric dendroidal spaces are a model for symmetric
\inftyoperad{}s (see \autoref{thm:cSegalmodelinfops} below). We briefly
retrace their main definitions in the world of \non{symmetric} operads. We will
use the resulting model category of complete Segal \emph{planar} dendroidal
spaces (or rather, its underlying \inftycat{}) as a model for
(\non{symmetric}) \inftyoperad{}s.

\begin{Cstr}~\cite[Sections 5 and 6]{CM13}
  \label{cstr:completeSegal}
  We build the simplicial model category $[\planarTreesop,\sSet]_\cSegal$ of
  \introduce{complete Segal dendroidal spaces} (also called
  \introduce{dendroidal Rezk model category}) as constructed by
  \nameCM{} in the symmetric case:\\

  Take the Reedy model structure\footnote{
    \nameCM{} actually use a
    generalized version of the Reedy model structure since the category
    $\symTrees$ of symmetric rooted trees is not a Reedy category (unlike
    $\planarTrees$, which is).
  }
  on the functor category $\dsSet\coloneqq
  [\planarTreesop,\sSet]$ and then
  Bousfield-localize~\cite[Proposition A.3.7.3]{Lurie2009} two times:
  \begin{enumerate}
  \item at the Segal core inclusions $\Segalcore T\lra \planarTrees[T]$ and
  \item at the maps $\planarTrees[T]\otimes \Jdendroidal\lra \planarTrees[T]$,
    where $\Jdendroidal$ is the dendroidal nerve of the category
    $\bullet\xra{\cong}\bullet$ with two objects and a single isomorphism
    between them.\qedhere
  \end{enumerate}
\end{Cstr}

The Reedy model category $[\planarTreesop,\sSet]_\Reedy$ has a canonical
simplicial enrichment~\cite[Theorem 10.3]{RV14} which is maintained by
the Bousfield localization processes~\cite[Proposition A.3.7.3]{Lurie2009}.
Therefore we can construct what we call the
\introduce{\inftycat{} of \inftyoperad{}s} as the simplicial nerve of
the fibrant-cofibrant objects:
\begin{equation}
  \infops\coloneqq \snerve{[\planarTreesop,\sSet]^\fibcofibdec_\cSegal}
\end{equation}

The name is justified by the following result.

\begin{Thm}~\cite[Corollary 6.8]{CM13}
  \label{thm:cSegalmodelinfops}
  The inclusion $\dSet\hra [\symTrees, \sSet]_\cSegal$ is a left Quillen
  equivalence between the model category of symmetric \inftyoperad{}s as
  defined by
  \nameCM{}~\cite{CM11}
  and the model category of complete Segal symmetric dendroidal spaces.
\end{Thm}

\begin{Def}\label{def:modelcatiSegal}
  We denote by $[\planarTreesop,\sSet]_\iSegal$
  the Bousfield localization of
  $[\planarTreesop,\sSet]_\cSegal$
  at the collapse maps
  \begin{equation}
    \planarTrees[\corolla n]\lra \planarTrees[T]
  \end{equation}
  for each $n$-ary tree $T$; we call it the \introduce{model category of invertible Segal dendroidal spaces}.\\
  We denote by
  \begin{equation}
    \invinfops\coloneqq
      \snerve{[\planarTreesop,\sSet]^\fibcofibdec_\iSegal}
  \end{equation}
  the corresponding \introduce{\inftycat{} of invertible \inftyoperad{}s}
\end{Def}
\begin{Rem}
  It is immediate from the characterization of Bousfield localization that
  $[\planarTreesop,\sSet]^\fibcofibdec_\iSegal$ is a full simplicial subcategory
  of $[\planarTreesop,\sSet]^\fibcofibdec_\cSegal$. Hence the
  \inftycat{} $\invinfops$ of invertible \inftyoperad{}s is a full
  subcategory of the \inftycat{} $\infops$ of (all) \inftyoperad{}s.
\end{Rem}

\begin{Lem}\label{lem:invertibleoperadsfullsubcat}
  The \inftycat{} $\invinfops$ of invertible \inftyoperad{}s is
  equivalent to the full subcategory of $\Fun(\planarTreesop, \Spaces)$
  consisting of those dendroidal spaces $\localX\colon \planarTreesop\to
  \Spaces$ which are invertible Segal and satisfy the following
  \introduce{completeness} condition:
  \begin{itemize}
  \item\label{item:inlem:complete} For each tree $T$,
    the map
    $\planarTrees[T]\otimes \Jdendroidal \to \planarTrees[T]$ from
    \autoref{cstr:completeSegal} induces an equivalence
    \begin{equation}
      \Hom(\planarTrees[T]\otimes \Jdendroidal,\localX)
      \xra{\simeq} \localX_T.
    \end{equation}
    in $\Spaces$.\qedhere
  \end{itemize}
\end{Lem}

To prove \autoref{lem:invertibleoperadsfullsubcat} we use the following result.
\begin{Prop}~\cite[Proposition 4.2.4.4.]{Lurie2009} %
  \label{prop:diagramsinmodelinftycat}
  \newcommand\localdiagcat{D}
  Let $\exmodelcat$ be a combinatorial simplicial model category, $\localdiagcat$ a
  small simplicial category and $S$ a simplicial set equipped with an
  equivalence $\catrealization S \xra{\simeq}\localdiagcat$.
  Then the induced map
  \begin{equation}
    \snerve{[\localdiagcat,\exmodelcat]^\fibcofibdec}
    \lra
    \Fun(S, \snerve{{\exmodelcat}^\fibcofibdec})
  \end{equation}
  is a categorical equivalence of simplicial sets.
\end{Prop}

\begin{Rem}
  \newcommand\localdiagcat{D}
  In \autoref{prop:diagramsinmodelinftycat} it does not matter whether we equip
  $[\localdiagcat,\exmodelcat]$ with the injective, projective or
  (if $\localdiagcat$ is a Reedy category)
  with the Reedy model structure, since they are all Quillen
  equivalent~\cite[Remark A.2.9.23]{Lurie2009}.
\end{Rem}

\begin{Prf}[of \autoref{lem:invertibleoperadsfullsubcat}]
  \newcommand\localdiagcat{D}
  We specialize \autoref{prop:diagramsinmodelinftycat} to $\exmodelcat \coloneqq
  \sSet$ and $\localdiagcat\coloneqq\planarTreesop$
  (seen as a discrete simplicial category);
  we put $S\coloneqq \nerve\planarTreesop = \snerve\planarTreesop$
  equipped with the adjunction counit
  $\catrealization{\snerve{\planarTreesop}}\xra{\simeq}\planarTrees$. We obtain
  an equivalence
  \begin{align}
    \label{equ:functorequivReedylevel}
    \snerve{[\planarTreesop,\sSet]_\Reedy^\fibcofibdec}\xra{\simeq}
    \Fun(\nerve\planarTreesop, \Spaces)
  \end{align}
  of \inftycats{}.
  Passing to Bousfield localizations replaces the simplicial category
  $[\planarTreesop,\sSet]^\fibcofibdec_\Reedy$
  by the full subcategory of the new fibrant-cofibrant objects.
  Therefore the equivalence \eqref{equ:functorequivReedylevel}
  restricts to an equivalence between
  $\invinfops\coloneqq\snerve{[\planarTreesop,\sSet]_\iSegal^\fibcofibdec}$
  and some full subcategory of
  $\Fun(\nerve\planarTreesop, \Spaces)$
  whose objects are determined
  by the fibrancy conditions in the three localization steps.
  Each of these steps corresponds precisely to one of the three conditions
  (invertibility, Segal, completeness)
  in \autoref{lem:invertibleoperadsfullsubcat}.
\end{Prf}

We will now see that the completeness condition in
\autoref{lem:invertibleoperadsfullsubcat} is redundant.

\begin{Lem}\label{lem:invertiblecompletenessisredundant}
  An invertible Segal dendroidal space is automatically complete.
\end{Lem}

\begin{Prf}
  A dendroidal Segal space
  $\localX\colon\planarTreesop\to\Spaces$ is complete if
  and only the underlying simplicial Segal space
  $\restr{\localX}{\Dop}\colon\Dop\subset\planarTreesop\to\Spaces$ (obtained by
  restricting to linear trees) is complete. If $\localX$ is invertible then
  $\restr{\localX}{\Dop}$ is constant, hence trivially complete.
\end{Prf}

\autoref{lem:invertiblecompletenessisredundant} motivates the name
\qquote{invertible Segal} (rather than \qquote{invertible complete Segal}) in
\autoref{def:modelcatiSegal} and completes the transition from
\autoref{cor:Segaldendrsimplspaces} to \autoref{cor:invinfopsSegalspaces}.

\begin{Rem}
\label{rem:main-for-Gamma}
The story of \autoref{sec:invertibleinfoperads}
can be retold, mutatis mutandis, in the world of
symmetric \inftyoperads{},
symmetric dentroidal spaces and $\Gamma$-spaces;
hence we obtain an equivalence between
the \inftycats{} of
\begin{itemize}
\item
  \twoSegal{} $\Gamma$-spaces and
\item
  invertible symmetric \inftyoperads{}.
  \qedhere
\end{itemize}
\end{Rem}

\begin{Rem}
  \autoref{expl:monoidal-grpd-cov-fibrant}
  and
  \autoref{rem:operads-monoidal-grpds} 
  have obvious analogs in the world
  of symmetric \inftyoperads{} and reduced Segal
  (\aka special) $\Gamma$-spaces.
\end{Rem}

\section{Variant: Cyclic operads and cyclic objects}

\label{sec:cyclicTrees}

We recall the definition of \ConnesCC{} $\Lambda$.
\begin{Def}~\cite{Connes1983}
  To each natural number $n\in \BN$ corresponds an object $\n\in\Lambda$ which
  we interpret as the unit circle $\sphere 1$ in the complex plane with $n+1$
  many equidistant marked points. The morphisms are homotopy classes of weakly
  monotone maps $\sphere 1\to \sphere 1$ of degree $1$ that send marked points
  to marked points.
\end{Def}
\begin{Rem}
  We fix the inclusion $\Delta\hra\Lambda$ which arranges the $n+1$ many
  elements of an object $\n\in\Delta$ as marked points on a circle. This
  inclusion is dense and faithful but not full.
\end{Rem}

We define the category $\cycTrees$ of \introduce{plane \emph{rootable} trees}.
In analogy to how $\planarTrees$ is a full subcategory of the category
$\Operads$ of operads, we define $\cycTrees$ as a full subcategory of the
category of cyclic operads
whose definition due to \twonames{Getzler}{Kapranov}\footnote{
  \twonames{Getzler}{Kapranov}
  introduced cyclic operads
  in their mono-colored and symmetric version.
}~\cite{GetzlerKapranov1995}
we now recall briefly.

\begin{Def}
  A \introduce{cyclic structure} on an operad $(\localO,O, \ocomp)$ consists of
  \begin{itemize}
	\item an involution $\ordual{(\blank)}\colon O\to O$ on colors (called
    \introduce{duality}) and
	\item a system of \introduce{rotation isomorphisms}
    \begin{equation}\localO(x_1,\ldots,x_n; y)\xra{\cong} \localO(\odual y, x_1, \ldots,
      x_{n-1}; \ordual{x_n})\end{equation} which is compatible with the composition of
    operations;
  \end{itemize}
  such that for each $n\in \BN$ the $(n+1)$-fold composition
  \begin{align*}
    \localO(x_1,\ldots,x_n; y)&\xra{\cong} \localO(\odual y, x_1, \ldots, x_{n-1}; \ordual{x_n})\xra{\cong} \localO(\ordual x_n, \ordual{y}, x_1, \ldots, x_{n-2}; \ordual{x_{n-1}})\\
                              &\xra{\cong}\cdots \xra{\cong}\localO(x_2,\ldots,x_n, \ordual y; \ordual{x_1})\xra{\cong}\localO(x_1,\ldots,x_n; y)
  \end{align*}
  of rotation isomorphisms is equal to the identity.

  A \introduce{cyclic operad} is an operad together with a cyclic structure. The
  cyclic operads are assembled into a category $\cycOperads$ where the morphisms
  are required to be compatible with the additional structure in the obvious
  way.
\end{Def}

\begin{Rem}
  We have an adjunction $\Operads\adjarrows\cycOperads$ where the right adjoint
  forgets the cyclic structure and the left adjoint adds a cyclic structure
  freely.
\end{Rem}

\begin{Def}A \introduce{plane rootable tree} consists of vertices and
  (unoriented) edges arranged in the plane, where an edge can connect two
  vertices or go to infinity in one or (in the case of the unique tree
  $\degtree$ with no vertices) both directions. We require our trees to have at
  least one external edge (this is what we mean by \qquote{rootable}). We think
  of each unoriented edge as a pair of anti-parallel arrows.
\end{Def}

\begin{Expl}
  A typical example of a plane rootable tree looks as follows:
  \[
    \begin{gathered}[b]
      \begin{tikzcd}
        &\,\doublear{ddr}&\,\doublear{rd}&\,\doublear{d}&\bullet\doublear{ld}\\
        \,\doublear{r}&\bullet\doublear{d}&&\bullet\doublear{d}&&\\
        \,\doublear{r}&\doublear{r}\doublear{ld}\bullet \doublear{r}&\bullet\doublear{ld}\doublear{r}\doublear{d}&\bullet\doublear{r}\doublear{d}&\bullet\doublear{r}&\,\\
        \,&\,&\bullet&\,&&
      \end{tikzcd}\\[-\dp\strutbox]
    \end{gathered}\qedhere
  \]
\end{Expl}

\newcommand\tsource[1]{s(#1)} \newcommand\ttarget[1]{t(#1)} We call an arrow a
\introduce{leaf} if comes from infinity and a \introduce{root} if it goes to
infinity. An arrow $a$ is called a \introduce{direct predecessor} of an arrow
$b$ (and $b$ is then a \introduce{direct successor} of $a$) if there is a vertex
which is both the target $\ttarget a$ of $a$ and the source $\tsource b$ of $b$.
We say $a$ is a \introduce{predecessor} of $b$ (or $b$ is a successor of $a$),
if $a$ is an iterated direct predecessor of $b$
(this includes the case $a=b$).
The \introduce{arity} of a tree (\resp a vertex) is
$n$, where $n+1$ is the number of arrows leaving (or, equivalently, entering)
the tree (\resp the vertex).

\begin{Rem}\label{rem:linearorderonpredleaves}
  For every arrow $b$ in a tree $T$, the set of predecessors of $b$ in $T$ forms
  a plane rooted tree (the root is $b$ itself). In particular there is a
  preferred linear order (clockwise along the boundary) on the set of those
  leaves $a$ of $T$ which are predecessors of $b$.
\end{Rem}

\begin{Cstr}
  Each plane tree $T$ gives rise to a cyclic operad (also denoted $T$) as
  follows:
  \begin{itemize}
  \item Each arrow is a color.
  \item Each pair $(v,a)$ consisting of an $n$-ary vertex $v\in T$ and an arrow
    $a$ starting in $v$ gives rise to an $n$-ary operation
    \begin{equation}v_a\colon( a_1,\dots, a_n) \lra a\end{equation} where the $a_i$'s are the direct
    predecessors of $a$ (hence $\ttarget {a_i}=v$) in clockwise order. All other
    operations are freely generated by these $v_a$'s.
  \item The involution on the colors exchanges the two anti-parallel arrows
    associated to a single edge.
  \item The rotation isomorphisms are given on generators by $v_a\mapsto
    v_{\ordual{a_n}}$.\qedhere
  \end{itemize}
\end{Cstr}

\begin{Def}
  We define the category $\cycTrees\subset\cycOperads$ of plane
  rootable trees to be the full subcategory spanned by the cyclic operads $T$
  constructed as above.
  A cyclic dendroidal object in an \inftycat{} $\localC$
  is a functor $\cycTreesop\to\localC$.
\end{Def}

\begin{Rem}
  Our category $\cycTrees$ is very close to the category of plane unrooted trees
  introduced by \twonames{Joyal}{Kock}~\cite{JK09};
  the only difference is that we require our trees to have at least one external
  edge. For instance, we do not allow the tree $\bullet$ which consists only of
  a single vertex, since this tree can not be interpreted as a cyclic operad in
  a meaningful way.
\end{Rem}

\begin{Rem}
  The free-cyclic-structure functor $\Operads\to \cycOperads$ induces an
  inclusion $\planarTrees\to\cycTrees$ which replaces each edge with two
  anti-parallel arrows and forgets the root.
\end{Rem}

\begin{Rem}
  The cyclic operad corresponding to the tree $\degtree$ (which has no vertices
  and exactly two mutually anti-parallel arrows) consists of two colors which
  are dual to each other and no \non{identity} operations. This cyclic operad
  $\degtree$ has an involution given by exchanging the two colors, i.e.\ the two
  arrows. A morphism $\degtree\to\localO$ to some cyclic operad $\localO$
  corresponds to a color of $\localO$; the involution on the colors of $\localO$
  is induced by the involution on $\degtree$.
\end{Rem}

\begin{Rem}
  It is easy to check that an operation in the cyclic operad $T\in\cycTrees$ is
  uniquely determined by its input and output colors. Hence a map $S\to T$
  between such operads is uniquely determined by the value at each arrow. Such a
  map would not, however, be determined by its values on unoriented edges; for
  instance, every unoriented edge $e$ of a tree $T$ gives rise to two different
  maps $\eta\to T$ in $\cycTrees$ corresponding to the two mutually dual colors
  described by $e$.

  If one were only interested in mono-colored cyclic operads or, more generally,
  cyclic operads with trivial duality (i.e.\ every color is self-dual), then it
  would be enough to consider unoriented edges. This point of view is taken by
  Hackney-Robertson-Yau~\cite{HRY19}.
\end{Rem}

\begin{Def}
  A map of plane rootable trees is called \introduce{boundary preserving} if it
  maps leaves to leaves and roots to roots. A \introduce{collapse map} in
  $\cycTrees$ is a boundary preserving map $\corolla{}\to T$ out of a
  corolla.
  A cyclic dendroidal object
  $\cycTreesop\to\localC$
  in some \inftycat{} $\localC$ is called \introduce{invertible} if
  it maps all collapse maps to equivalences in $\localC$.
\end{Def}

As the notation suggests, the category $\cycTrees$ of plane rootable trees has a
close relationship to the cyclic category: the latter is a localization of the
former as we will see next.

\begin{Cstr}[Covariant description of $\treerootcyc$]
  Analogously to the case of plane rooted trees, a plane rootable tree
  partitions the plane into \roughly{areas} which are arranged clockwise around
  a circle. This assignment is a functor $\treerootcyc\colon\cycTrees\to
  \Lambda$ which extends the functor $\treerootpl\colon\planarTrees \to
  \Delta$.\end{Cstr}

\begin{Cstr}[Contravariant description of $\treerootcyc$]
  Using the self-duality $\Lambda\cong\Lambda^\op$ (which interchanges marked
  points and intervals on a circle) we can define the functor
  $L\colon\cycTrees\to\Lambda^\op$ instead:

  A tree $T$ gets mapped to its set of leaves which are naturally arranged
  around a circle. The image of a morphism $\alpha\colon S\to T$ sends each leaf
  $a$ of $T$ to the unique leaf $b$ of $S$ such that $\alpha(b)$ is a successor
  of $a$. This assignment does not yet uniquely determine $L\alpha$ as a
  morphism in $\Lambda$; we still need to specify a linear order on the
  pre-images $(L\alpha)^\inv(b)$ (for every leaf $b$ of $S$) but this is taken
  care of by \autoref{rem:linearorderonpredleaves}.
\end{Cstr}

\begin{Rem}
  By combining the ideas from \autoref{sec:cyclicTrees} and
  \autoref{sec:symmetrictrees} we can construct a category of (\non{plane})
  rootable trees as a full subcategory of cyclic symmetric operads\footnote{
    Such
    operads have both a cyclic and a symmetric structure which are compatible
    when regarding the symmetric group $\Sn$ and the cyclic group
    $\cygrp{(n+1)}$ as a subgroup of $\SG {n+1}$.
  }.
  The corresponding functor
  $\treerootabstract\colon\abstractTrees\to\finsetneop$ maps a tree to its
  \non{empty} set of leaves (i.e.\ incoming arrows).
\end{Rem}

\autoref{prop:situationC} still holds for
$\treeroots\in\set{\treerootcyc, \treerootabstract}$
with essentially the same proof,
hence \autoref{lem:quillentheoremC} yields the following
cyclic version of \autoref{thm:localizationW}:

\begin{dThm}\label{thm:localizationWL}
  The functor
  $\treerootcyc\colon\cycTrees\to\Lambda$
  (\resp $\treerootabstract\colon\abstractTrees\to\finsetneop$)
  exhibits $\Lambda$
  (\resp $\finsetneop$)
  as an \inftycategorical{} localization of
  $\cycTrees$
  (\resp $\abstractTrees$)
  at the set of collapse maps.
\end{dThm}

\begin{Cor}
  The classifying space of $\cycTrees$ is $\classB{\sphere 1}$.
\end{Cor}
\begin{Prf}
  Follows immediately from \autoref{thm:localizationWL}
  because the classifying space of the cyclic category $\Lambda$
  is known to be $\classB{\sphere 1}$~\cite[Theorem 10]{Connes1983}.
\end{Prf}

\begin{Rem}
  \label{rem:careful-with-Lambda}
  Analogously to \autoref{cor:Segaldendrsimplspaces} 
  one can show that the functor $\treerootcyc$
  induces an equivalence between the \inftycats{} of
  \begin{itemize}
  \item
    \twoSegal{} cyclic objects and
  \item
    invertible cyclic  Segal dendroidal objects
  \end{itemize}
  in any \inftycat{} $\localC$,
  where \twoSegal{}/Segal are defined either as the obvious analogs of
  \autoref{def:Segalgrafting} and \autoref{def:unital2Segal}
  or, alternatively, by referring to the underlying simplicial/dendroidal object.

  Unfortunately, there is currently no result in the literature
  exhibiting (complete) Segal cyclic dendroidal spaces
  as a model for cyclic \inftyoperads{}.
  One promising approach to resolve this issue is proposed by \nameDCH{}
  who construct~\cite[Theorem 6.5]{DH18}
  a Dwyer--Kan type model structure on the category of
  simplicially enriched cyclic operads\footnote{
    \nameDCH{} call \buzzword{non-$\Sigma$ positive cyclic operads}
    what we simply call cyclic operads.}
  and conjecture\cite[Remark 6.9]{DH18} that it should be Quillen equivalent
  to a \roughly{complete Segal space}-type
  model structure on cyclic dendroidal simplicial sets
  lifted from the complete Segal model structure of \nameCM{}.
  Conditional on their conjecture,
  we can then say that \twoSegal{} cyclic spaces
  are equivalent to
  invertible cyclic \inftyoperads{}.
\end{Rem}


\newcommand{\etalchar}[1]{$^{#1}$}
\providecommand{\bysame}{\leavevmode\hbox to3em{\hrulefill}\thinspace}
\providecommand{\MR}{\relax\ifhmode\unskip\space\fi MR }
\providecommand{\MRhref}[2]{%
  \href{http://www.ams.org/mathscinet-getitem?mr=#1}{#2}
}
\providecommand{\href}[2]{#2}
\begin{thebibliography}{GCKT18b}

\bibitem[ACM19]{ACM2019}
Dimitri Ara, Denis-Charles Cisinski, and Ieke Moerdijk, \emph{The dendroidal
  category is a test category}, Math. Proc. Cambridge Philos. Soc. \textbf{167}
  (2019), no.~1, 107--121.

\bibitem[BM17]{BM17}
Pedro {Boavida de Brito} and Ieke {Moerdijk}, \emph{Dendroidal spaces,
  {$\Gamma$}-spaces and the special {Barratt-Priddy-Quillen} theorem}.

\bibitem[BOO{\etalchar{+}}18]{BOORS2018}
Julia~E. {Bergner}, Ang{\'e}lica~M. {Osorno}, Viktoriya {Ozornova}, Martina
  {Rovelli}, and Claudia~I. {Scheimbauer}, \emph{{2-Segal objects and the
  Waldhausen construction}}.

\bibitem[BV73]{BV73}
J.~M. Boardman and R.~M. Vogt, \emph{Homotopy invariant algebraic structures on
  topological spaces}, Lecture Notes in Mathematics, Vol. 347, Springer-Verlag,
  Berlin-New York, 1973.

\bibitem[CM11]{CM11}
Denis-Charles Cisinski and Ieke Moerdijk, \emph{Dendroidal sets as models for
  homotopy operads}, J. Topol. \textbf{4} (2011), no.~2, 257--299.

\bibitem[CM13]{CM13}
\bysame, \emph{Dendroidal {Segal} spaces and {$\infty$}-operads}, J. Topol.
  \textbf{6} (2013), no.~3, 675--704.

\bibitem[Con83]{Connes1983}
Alain Connes, \emph{Cohomologie cyclique et foncteurs {${\rm Ext}^n$}}, C. R.
  Acad. Sci. Paris S\'{e}r. I Math. \textbf{296} (1983), no.~23, 953--958.

\bibitem[DH18]{DH18}
Gabriel~C. {Drummond-Cole} and Philip {Hackney}, \emph{{Dwyer--Kan homotopy
  theory for cyclic operads}}.

\bibitem[DK12]{DyckerhoffKapranov2012}
Tobias {Dyckerhoff} and Mikhail {Kapranov}, \emph{{Higher {S}egal spaces I}}.

\bibitem[DK18]{DyckerhoffKapranov2018}
\bysame, \emph{Triangulated surfaces in triangulated categories}, J. Eur. Math.
  Soc. (JEMS) \textbf{20} (2018), no.~6, 1473--1524.

\bibitem[Dyc18]{Dyckerhoff2018}
Tobias Dyckerhoff, \emph{Higher categorical aspects of {H}all algebras},
  Building bridges between algebra and topology, Adv. Courses Math. CRM
  Barcelona, Birkh\"{a}user/Springer, Cham, 2018, pp.~1--61.

\bibitem[FGK{\etalchar{+}}19]{FGKUPW2019}
Matthew {Feller}, Richard {Garner}, Joachim {Kock}, May~U. {Proulx}, and Mark
  {Weber}, \emph{{Every 2-Segal space is unital}}.

\bibitem[GCKT18a]{GCKT2018a}
Imma G\'{a}lvez-Carrillo, Joachim Kock, and Andrew Tonks, \emph{Decomposition
  spaces, incidence algebras and {M}\"{o}bius inversion {I}: {B}asic theory},
  Adv. Math. \textbf{331} (2018), 952--1015.

\bibitem[GCKT18b]{GCKT2018b}
\bysame, \emph{Decomposition spaces, incidence algebras and {M}\"{o}bius
  inversion {II}: {C}ompleteness, length filtration, and finiteness}, Adv.
  Math. \textbf{333} (2018), 1242--1292.

\bibitem[GCKT18c]{GCKT2018c}
\bysame, \emph{Decomposition spaces, incidence algebras and {M}\"{o}bius
  inversion {III}: {T}he decomposition space of {M}\"{o}bius intervals}, Adv.
  Math. \textbf{334} (2018), 544--584.

\bibitem[GK95]{GetzlerKapranov1995}
E.~Getzler and M.~M. Kapranov, \emph{Cyclic operads and cyclic homology},
  Geometry, topology, \& physics, Conf. Proc. Lecture Notes Geom. Topology, IV,
  Int. Press, Cambridge, MA, 1995, pp.~167--201.

\bibitem[HRY19]{HRY19}
Philip Hackney, Marcy Robertson, and Donald Yau, \emph{Higher cyclic operads},
  Algebr. Geom. Topol. \textbf{19} (2019), no.~2, 863--940.

\bibitem[JK09]{JK09}
Andr{\'e} {Joyal} and Joachim {Kock}, \emph{{Feynman} graphs, and nerve theorem
  for compact symmetric multicategories (extended abstract)}.

\bibitem[Joy02]{Joyal2002}
A.~Joyal, \emph{Quasi-categories and {K}an complexes}, J. Pure Appl. Algebra
  \textbf{175} (2002), no.~1-3, 207--222, Special volume celebrating the 70th
  birthday of Professor Max Kelly.

\bibitem[Joy08]{Joyal2008}
Andr{\'e} Joyal, \emph{Notes on quasi-categories}, Lecture Notes, 2008.

\bibitem[KS11]{KS11}
Maxim Kontsevich and Yan Soibelman, \emph{Cohomological {H}all algebra,
  exponential {H}odge structures and motivic {D}onaldson-{T}homas invariants},
  Commun. Number Theory Phys. \textbf{5} (2011), no.~2, 231--352.

\bibitem[Lur09]{Lurie2009}
Jacob Lurie, \emph{Higher topos theory}, Annals of Mathematics Studies, vol.
  170, Princeton University Press, Princeton, NJ, 2009.

\bibitem[Lur17]{Lurie2017}
\bysame, \emph{Higher algebra}, September 2017.

\bibitem[May72]{May72}
J.~P. May, \emph{The geometry of iterated loop spaces}, Springer-Verlag,
  Berlin-New York, 1972, Lectures Notes in Mathematics, Vol. 271.

\bibitem[MW07]{MW07}
Ieke Moerdijk and Ittay Weiss, \emph{Dendroidal sets}, Algebr. Geom. Topol.
  \textbf{7} (2007), 1441--1470.

\bibitem[Qui67]{Quillen1967}
Daniel~G. Quillen, \emph{Homotopical algebra}, Lecture Notes in Mathematics,
  No. 43, Springer-Verlag, Berlin-New York, 1967.

\bibitem[Rez01]{Rezk2001}
Charles Rezk, \emph{A model for the homotopy theory of homotopy theory}, Trans.
  Amer. Math. Soc. \textbf{353} (2001), no.~3, 973--1007.

\bibitem[RV14]{RV14}
Emily Riehl and Dominic Verity, \emph{The theory and practice of {R}eedy
  categories}, Theory Appl. Categ. \textbf{29} (2014), 256--301.

\bibitem[Seg74]{Segal1974}
Graeme Segal, \emph{Categories and cohomology theories}, Topology \textbf{13}
  (1974), 293--312.

\bibitem[{Ste}19]{Stern2019}
Walker~H. {Stern}, \emph{{2-Segal objects and algebras in spans}}.

\bibitem[To{\"{e}}06]{Toen2006}
Bertrand To{\"{e}}n, \emph{Derived {H}all algebras}, Duke Math. J. \textbf{135}
  (2006), no.~3, 587--615.

\bibitem[Wal85]{Waldhausen1985}
Friedhelm Waldhausen, \emph{Algebraic {$K$}-theory of spaces}, Algebraic and
  geometric topology ({N}ew {B}runswick, {N}.{J}., 1983), Lecture Notes in
  Math., vol. 1126, Springer, Berlin, 1985, pp.~318--419.

\end{thebibliography}
\end{document}